\documentclass[11pt]{article}
\usepackage{amsmath,amsthm}
\usepackage{fontspec}
\usepackage{unicode-math}
\usepackage{graphicx}
\usepackage[margin=1.05in]{geometry}
\usepackage[hidelinks,pdfversion=1.7]{hyperref}
\newtheorem{theorem}{Theorem}[section]
\newtheorem{proposition}[theorem]{Proposition}
\newtheorem{lemma}[theorem]{Lemma}
\newtheorem{corollary}[theorem]{Corollary}

\newcommand{\diag}{\operatorname{diag}}
\newcommand{\wind}{\operatorname{wind}}

\title{Critical spectral weights on triangulated disks}

\author{Enqiang Zhu\\ Institute of Computing Science and Technology,\\ Guangzhou University, Guangzhou 510006, China}

\begin{document}
\maketitle

\begin{abstract}
We solve the common-positive-weight support problem at the local spectral
threshold $1/2$ for relative triangulated disks with interior degrees at
most five and an empty or connected degree-five core. We characterize
feasibility by a finite decomposition into three elementary pieces and
a condition on the original boundary marking. The proof passes from
spectral equality to a positive conjugate rectangle, then uses geometric
reductions to obtain a serial presentation. The empty-face inequality
holds for some common positive weights exactly when such a presentation
has only old join endpoints. Whenever this condition holds, our
construction yields weights with a positive definite empty-face matrix.
Positive Schur responses establish sufficiency, and a uniform
face-counting obstruction proves necessity. This gives a finite
combinatorial criterion for critical spectral feasibility together
with a constructive choice of common weights.
\end{abstract}

\section{The support problem and the classification}
We determine when a relative triangulated disk admits one positive
triangle weighting that satisfies all vertex and empty-face spectral
conditions at the critical threshold $1/2$. For disks with interior
degrees at most five and an empty or connected degree-five core, we
give a finite combinatorial classification and construct the common
weights whenever the criterion is satisfied.

The central requirement is compatibility: the same facet weights
must induce all the link matrices. This places the problem within
the study of spectral supports motivated by Vinzant's Question~5.15,
which asks which simplicial complexes support local spectral expansion
beyond the log-concave setting \cite{Vinzant}. The classification below
identifies the geometric structure and boundary marking that enforce
this compatibility for the relative disk problem.

The method has two parts. First, the equality equations turn the
simultaneous vertex inequalities into a positive conjugate rectangle.
Removing its independent four-valent stars leaves a constrained planar
core. Four coordinate extrema bound the complexity of every core block;
geometric exclusions and two legal inverse flips then exhaust the
connected cores. Second, the reverse construction has one fixed axis.
At interfaces with old endpoints, scalar Schur responses construct
strict common weights. At the other interfaces, a uniform count of
face contributions gives an obstruction. The resulting support criterion is a finite test on the triangulation
and its original boundary marking.

The local-to-global viewpoint begins with Garland's spectral method
\cite{Garland}, its development for polyhedral complexes
\cite{BallmannSwiatkowski}, and weighted simplicial formulations
\cite{Oppenheim}. For a pure weighted simplicial complex, conditioning
on a face gives a weighted link graph. Local spectral expansion asks
for bounds on the nonconstant spectrum of these graphs, with all link
weights induced by the same facet measure. Such bounds can control
global Laplacians and higher-order random walks
\cite{OppenheimDescent,KaufmanOppenheim}.
Compatibility of the link weights is therefore built into the facet
measure.

Log-concave polynomials provide an important class of common weights.
Anari, Oveis Gharan and Vinzant developed the polynomial framework
in \cite{AOV}. Anari, Liu, Oveis Gharan and Vinzant \cite{ALOV} relate
localized spectra to derivatives of generating polynomials and prove
rapid mixing for matroid bases. Lorentzian polynomials provide a
related structural theory of matroidal supports \cite{BrandenHuh}.
These results explain the log-concave regime behind the question in
\cite{Vinzant}. Expansion of complexes also has applications quite
different from support classification, including the construction of
good locally testable codes \cite{DinurEtAl}.

There is also a geometric comparison. Positive barycentric embeddings
go back to Tutte \cite{Tutte}; degree methods for mappings of meshes
with boundary are developed in \cite{Lipman}. Electrical conjugacy
and rectangle tilings form a separate classical source of geometric
coordinates, with recent convergence results in \cite{BP}. In another
discrete-geometric setting, Ge and Lin \cite{GeLin} obtain degree
criteria for hyperbolic circle packings using combinatorial Ricci
flow. In the present construction, spectral equality produces
coordinates with an indefinite orthogonality relation. Indefinite forms
also occur in Moutard-net theory \cite{BS}. Schur complementation and
network elimination supply the interface algebra \cite{JLLDS}.

Our contribution is the passage from spectral equality to a complete
marked support classification. Positive conjugate geometry forces
the permitted reductions, their reversal gives the finite presentations,
and the original interface marking decides the empty-face condition.
Together, these steps turn a simultaneous spectral feasibility problem
into a combinatorial criterion with a constructive positive-weight
witness.

\subsection*{Relative operators}
Let $C$ be a finite simplicial triangulation of an oriented closed disk.
Let $L=L_0\sqcup L_1$ consist of two vertex-disjoint nondegenerate
boundary paths. A vertex or edge is \emph{old} if it belongs to $L$,
and \emph{new} otherwise. Put $N=C^0\smallsetminus L^0$. For positive
triangle weights $t=(t_T)$ define
\[
 d_{uv}=\sum_{T\supset uv}t_T,\qquad m_v=\sum_{T\ni v}t_T.
\]
The symmetric normalized adjacency of the link of $v$ is
\[
 S_v(u,w)=\frac{t_{vuw}}{\sqrt{d_{vu}d_{vw}}},
\]
with missing entries and the diagonal zero. At an old vertex let $Q_v$
be its principal submatrix on coordinates $u$ with $vu\notin L$.
These submatrices retain the degrees from the full link. For the empty face
retain the original new set and put
\[
 Q_\varnothing(u,v)=\frac{d_{uv}}{2\sqrt{m_um_v}}
 \qquad(u,v\in N,\ u\ne v).
\]
Empty matrices impose no constraint. With eigenvalues in decreasing
order, write
\[
 \begin{aligned}
 \alpha(C,L;t)=\max\{0,&\max_{v\in N}\lambda_2(S_v),
 \max_{v\in L^0}\lambda_{\max}(Q_v),
 \lambda_{\max}(Q_\varnothing)\},\\
 b(C,L)&=\inf_{t>0}\alpha(C,L;t).
 \end{aligned}
\]
The \emph{vertex conditions} $V$ say that the two vertex maxima are at
most $1/2$. Complete criticality also imposes the empty-face inequality,
equivalently
\begin{equation}\label{eq:empty}
 H_N:=\diag(m_v)_{v\in N}-(d_{uv})_{u,v\in N}\succeq0,
\end{equation}
where the adjacency has zero diagonal. The boundary killing is encoded
by these principal submatrices with their full-link normalizations.

\subsection*{The three pieces}
\begin{figure}[htbp]
\centering
\includegraphics[width=\linewidth]{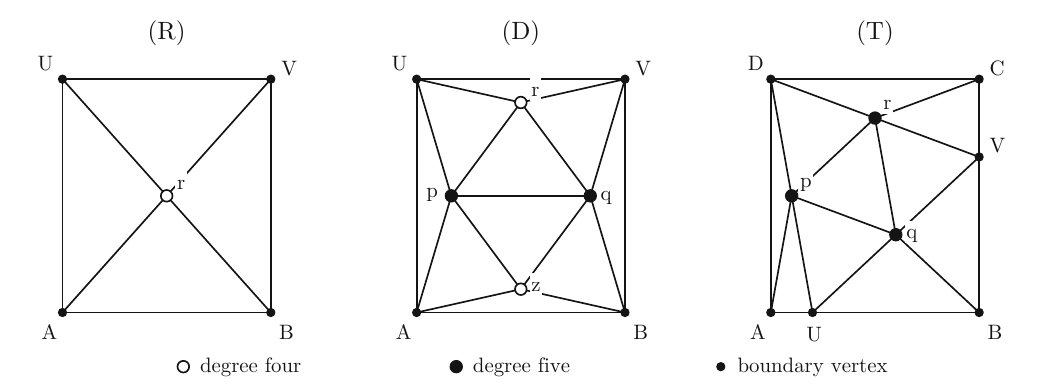}
\caption{The three pieces, with exactly the edges of the face lists.
Open interior circles have degree four and filled interior circles degree
five; small dots are boundary vertices. These drawings use positive
$E=0$ charts, independently scaled in the two axes. In T, $U,V$ are
strict flat boundary vertices, not additional corners.}
\label{fig:modules}
\end{figure}

The following disks have named boundary vertices; their names specify
combinatorics, not metric coordinates.
\begin{enumerate}
\item[(R)] A quadrilateral $A,B,V,U$ coned from an interior vertex $r$:
\[
 ABr,\ BVr,\ VUr,\ UAr.
\]
\item[(D)] A quadrilateral $A,B,V,U$ with interior vertices $p,q,z,r$:
\[
 \begin{gathered}
 ABz,\ BVq,\ VUr,\ UAp,\ Bzq,\ Vqr,\ Urp,\ Apz,\ zqp,\ pqr.
 \end{gathered}
\]
Here $z,r$ have degree four and $p,q$ have degree five.
\item[(T)] A six-cycle $A,U,B,V,C,D$ with interior vertices $p,q,r$:
\[
 DAp,\ AUp,\ UBq,\ BVq,\ VCr,\ CDr,\ pUq,\ qVr,\ rDp,\ pqr.
\]
Its four sides are $AUB,BVC,CD,DA$; its three interior degrees are five.
\end{enumerate}
Figure~\ref{fig:modules} gives geometric representatives.

An \emph{admissible presentation} is one of the following constructions.
Take a nonempty finite word of R and D pieces and glue successive pieces
along opposite full end edges, without any other identifications.
Either opposite side pair may be the ends of D. Designated joining edges
may then be flipped to the edge between their two incident interior
apices, provided both apices have degree four before the flip and no
apex is used twice. Alternatively, take T and attach a possibly empty
unflipped R strip on $DA$ or on $CD$, but not both. The old paths are
one pair of opposite complete outer side paths. The original endpoints
of every join remain marked even when its edge is flipped away.

Let $G_5$ be the graph induced by the original interior degree-five
vertices, including isolated vertices. Define $\mathcal P(C,L)$ to mean
that $(C,L)$ has an admissible presentation in which every original
join endpoint belongs to $L^0$. This is a finite combinatorial predicate:
the number of pieces is bounded by the number of triangles, and the
construction involves only identifications, markings and specified flips.

\begin{theorem}[Connected-core support classification]\label{thm:main}
Let $(C,L)$ be as above. Suppose all original interior degrees are at
most five and $G_5$ is empty or connected.
\begin{enumerate}
\item Common positive weights satisfying $V$ exist if and only if
$(C,L)$ has an admissible presentation.
\item The following are equivalent:
\[
 \begin{aligned}
 &\mathcal P(C,L);\\
 &\exists\,t_T>0:\quad V\ \text{and}\ H_N\succeq0;\\
 &\exists\,t_T>0:\quad V\ \text{and}\ H_N\succ0.
 \end{aligned}
\]
\end{enumerate}
All matrices in each existence statement come from the same weights.
One may normalize $\sum_Tt_T=1$.
\end{theorem}
For each support satisfying the criterion, the construction supplies a
common weighting with $H_N\succ0$. A one-piece presentation has no joins,
so its interface condition is vacuous. For two R pieces with their join
flipped, denoted F below, the original join endpoints enter the same
marking test.

\section{Criticality and conjugate rectangles}
For a vertex function \(x\), with \(x=0\) on \(L_0\) and \(x=1\)
on \(L_1\), consider the midpoint equations
\begin{equation}\label{eq:M}
 \sum_{T=uvw}t_T\left(x_w-\frac{x_u+x_v}{2}\right)=0
 \qquad(uv\notin L).
\end{equation}

\begin{proposition}\label{prop:rectangle}
The conditions \(V\) imply \eqref{eq:M} for the unique weighted
Dirichlet potential \(x\). Conversely, any positive solution of
\eqref{eq:M} gives a piecewise affine orientation-preserving
homeomorphism
\[
 F=(x,y):C\longrightarrow[0,1]\times[0,H],\qquad H>0,
\]
with old paths vertical and complementary boundary paths horizontal.
On each abstract unit equilateral face, with \(J\) denoting
counterclockwise rotation by \(\pi/2\),
\begin{equation}\label{eq:conjugacy}
 \nabla y=t_TJ\nabla x .
\end{equation}
All old and new boundary vertex conditions then hold. If every
interior vertex has degree at most five, \(V\) is equivalent to
\eqref{eq:M}. For every positive weighting,
\begin{equation}\label{eq:lower}
 \alpha(C,L;t)\ge\tfrac12.
\end{equation}
\end{proposition}

\begin{proof}
Minimize \(\sum_{uv}d_{uv}(x_u-x_v)^2\) with the prescribed old
values. Every component of new vertices meets the old set, so the
minimizer is unique. Put \(f_v(u)=x_u-x_v\). At a new vertex its
\(d_{vu}\)-weighted mean is zero; at an old vertex its old-edge
coordinates vanish. Write \(A_v,D_v\) for the link adjacency and
degree matrices. Direct counting on each triangle gives
\begin{equation}\label{eq:garland}
 \sum_v f_v^{\mathsf T}A_v f_v
       =\tfrac12\sum_v f_v^{\mathsf T}D_v f_v .
\end{equation}
The vertex inequalities bound each summand on the left by its
corresponding summand on the right. Equality in a self-adjoint
quadratic bound therefore forces
\(A_vf_v=\frac12D_vf_v\) on every new-edge coordinate, which is
\eqref{eq:M}. The same identity rules out strict vertex inequalities
below \(1/2\), since \(x\) is not constant. This proves
\eqref{eq:lower} without using the empty face.

Now assume only \eqref{eq:M}. Summing the link eigen-equations at a
new vertex proves its Dirichlet harmonic equation, so \(0\le x\le1\).
A constant triangle would force its neighbor across each interior
edge to be constant by \eqref{eq:M}. Connectedness of the triangle
dual would then contradict the two old values. There are no constant
triangles. An extreme value \(0\) or \(1\) at a new vertex, by
harmonicity, would make its entire star constant. Thus
\(0<x_v<1\) at every new vertex.

On an oriented face \(uvw\) define
\begin{equation}\label{eq:dy}
 dy_{uv}=\frac{t_T}{\sqrt3}(x_u+x_v-2x_w).
\end{equation}
The face sum is zero; opposite orientations and \eqref{eq:M} make
the definitions agree on shared edges. Simple connectivity integrates
this cochain to \(y\). An equilateral-face calculation gives
\eqref{eq:conjugacy}, hence
\(\det DF=t_T|\nabla x|^2>0\).

In the induced boundary orientation \(dy<0\) on \(L_0\), \(dy>0\)
on \(L_1\), and \(dy=0\) on free edges. Consequently the boundary
lies on a rectangle of positive height, though free paths might
initially appear to backtrack. Their endpoint-fixed homotopy to the
straight horizontal segments shows boundary winding number one
inside the rectangle and zero outside. For a point off all image
edges, cancellation of interior edges gives
\[
 \wind(F(\partial C),p)
   =\sum_T\wind(F(\partial T),p)
   =\#\{T:p\in F(T)^\circ\}.
\]
Every image triangle is therefore inside the rectangle. A backwards
edge on its lower or upper side would put the positively oriented
incident triangle outside; hence both free paths are strictly
monotone. The boundary is a homeomorphism, and the same degree
identity, together with positive local triangle sectors, gives the
global PL homeomorphism.

At an old vertex the retained \(f_v\) has one strict sign.
Indeed a same-color old chord would make a constant face by
\eqref{eq:M}; new neighbors have strict intermediate values. The
retained link is a connected path, so its positive eigenfunction
gives Perron eigenvalue \(1/2\). At a new boundary vertex the embedded
half-star is cut once by the vertical line through its image.
Thus the link eigenfunction has exactly one sign change after zero
coordinates are removed. To justify the spectral consequence
including zeros, set \(K=\frac12D_v-A_v\). Zeros of \(f_v\) are
isolated and cannot be endpoints, by the positive path recurrence.
Their diagonal block in \(K\) is positive; its Schur complement
joins the two neighbors of a zero with positive weight
\(2pq/(p+q)\). Congruence by the remaining nonzero \(f_v\) turns this
Schur complement into a signed path Laplacian, with conductance
\(h_{ij}f_v(i)f_v(j)\). Its negative index equals the number of sign
changes, since path increments are independent modulo constants.
This number is one. Hence the original full link has at most one
eigenvalue greater than \(1/2\), as required.

At an interior degree-four vertex the weighted cycle spectrum is
\(1,\eta,-\eta,-1\). Its known eigenvalue \(1/2\) must be \(\eta\).
For a positive weighted five-cycle, all eigenvalues exceed \(-1\).
If its third eigenvalue were at least \(1/2\), the largest three
would sum to at least two and the last two to more than minus two,
contradicting zero trace. Thus its known half eigenvalue is the second.
Degrees below four are excluded by the same-sign neighbor argument
in Theorem~\ref{thm:fourcore}.
This completes the asserted equivalence.
\end{proof}

\begin{corollary}\label{cor:degree}
In a positive midpoint rectangle, a boundary vertex other than
one of the four old-path endpoints has degree at least four.
Each of those four endpoints has degree at least three.
\end{corollary}
\begin{proof}
At an old vertex the nonzero half-eigenfunction requires at least
two retained coordinates. There are two deleted coordinates in an
old path's interior and one at its endpoints. A new boundary link
of size two or three has spectrum \(1,-1\) or \(1,0,-1\),
respectively, and cannot have eigenvalue \(1/2\).
\end{proof}

\section{The four-valent layer}
For \(z=(z_x,z_y)\), put \(q(z)=z_xz_y\), and write
\(B(u,v)=u_xv_y+u_yv_x\). Set
\[
 E(a,b,c)=q(a-b)+q(b-c)+q(c-a).
\]
On a unit equilateral reference face, orthogonality of the two
coordinate gradients is exactly \(E=0\). Conversely, a positive
oriented triangle with \(E=0\) recovers its positive weight by
\begin{equation}\label{eq:weight}
 t=\frac{\sqrt3}{2}\frac{\det(p,r)}{p_x^2+r_x^2-p_xr_x},
 \qquad p=b-a,\quad r=c-a.
\end{equation}
The denominator is positive because positive area excludes
\(p_x=r_x=0\). Direct calculation gives \(\nabla y=tJ\nabla x\).
Consequently an embedded rectangle triangulation with positive
\(E=0\) faces supplies common midpoint weights, by agreement of
the edge increments in \eqref{eq:dy}.
\subsection*{The explicit four-sided rule}
For cyclic \(a,b,c,d\), set
\[
 U=a-c,\quad V=b-d,\quad D=\det(U,V),\quad
 m=(a+c)/2,\quad n=(b+d)/2,\quad \delta=m-n.
\]
When \(D\ne0\), write \(\delta=\mu U+\nu V\).

\begin{corollary}\label{cor:quad}
These four points bound a positive embedded \(E=0\) fan if and
only if
\begin{equation}\label{eq:quad}
 \begin{gathered}
 D>0,\qquad B(U,V)=0,\qquad q(U)+q(V)+12q(\delta)=0,\\
 (1+6|\mu|)(1+6|\nu|)<4 .
 \end{gathered}
\end{equation}
Its unique center is \(z=m+\mu U-2\nu V\), and its boundary
quadrilateral is strictly convex.
\end{corollary}
\begin{proof}
Differences of adjacent face equations give
\(B(U,V)=0\),
\(B(U,z)=B(U,2m-n)\), and \(B(V,z)=B(V,2n-m)\).
They force the displayed center. At that center each face
expression equals
\(\bigl(q(U)+q(V)+12q(\delta)\bigr)/2\).
The four consecutive determinants divided by \(D\) are
\[
 \begin{aligned}
 \gamma_1&=\tfrac14+\tfrac{\nu-\mu}{2}+3\mu\nu,&
 \gamma_2&=\tfrac14+\tfrac{\mu+\nu}{2}-3\mu\nu,\\
 \gamma_3&=\tfrac14+\tfrac{\mu-\nu}{2}+3\mu\nu,&
 \gamma_4&=\tfrac14-\tfrac{\mu+\nu}{2}-3\mu\nu .
 \end{aligned}
\]
They sum to one. The opposite pairs are positive exactly when
\(1+12\mu\nu>2|\mu-\nu|\) and
\(1-12\mu\nu>2|\mu+\nu|\). According to the sign of \(\mu\nu\),
the stronger of these is
\(2(|\mu|+|\nu|)+12|\mu\nu|<1\), the last condition in
\eqref{eq:quad}. Four positive successive angles have total
strictly less than \(4\pi\), so their positive integer winding
is one. Finally the paired sums of the \(\gamma_i\) give
\(|\mu|,|\nu|<1/2\). Both diagonals meet strictly inside their
segments, proving strict convexity.
\end{proof}

\subsection*{The paired-star normal form}
Consider the six-triangle disk with outer boundary \((u,c,v,d)\),
interior vertices \(a,b\), and faces
\(auc,acv,avb,abu,bvd,bdu\). Its four boundary vertices are distinct.
This is the union of two degree-four stars with a shared edge.

\begin{proposition}\label{prop:pair}
An ordered simple boundary quadrilateral admits a positive embedded
\(E=0\) realization of this disk only if it is a parallelogram.
Set \(o=(u+v)/2=(c+d)/2\), \(U=u-v\), \(C=c-o\).
Necessarily \(q(U)\ne0\). Define
\[
 h=-\frac{B(U,C)}{2q(U)},\qquad T=C+hU,\qquad
 D=\det(U,T).
\]
Then \(D>0\). Also \(B(U,T)=0\) and \(q(U)\ne0\) imply
\(q(T)\ne0\), so set
\(\Delta=-(1+12h^2)q(U)/q(T)-3\).
Every positive realization satisfies
\begin{equation}\label{eq:pairdomain}
 D>0,\qquad |h|<\tfrac16,\qquad
 0\le\Delta<\left(\frac{1-6|h|}{1+2|h|}\right)^2.
\end{equation}
Its two centers have the form
\begin{equation}\label{eq:paircenters}
 a=o+hU+\frac{1+\tau}{2}T,\qquad
 b=o-hU+\frac{\tau-1}{2}T,\qquad \tau^2=\Delta.
\end{equation}

\end{proposition}
\begin{proof}
Translate \((u+v)/2\) to zero. Each positive four-star has strictly
convex boundary by Corollary~\ref{cor:quad}. Its diagonals \(U,c-b\)
are independent and \(B\)-orthogonal, so
\(4q(U)q(c-b)=-\det(U,c-b)^2\ne0\).
The shared-face difference gives \(B(U,a+b)=0\); the two
outer differences give \(B(U,a+c)=B(U,b+d)=0\).
In a basis \(U,T_0\) with \(B(U,T_0)=0\), write
\[
 a=hU+kT_0,\quad b=-hU+lT_0,\quad
 c=-hU+mT_0,\quad d=hU+nT_0.
\]
The four-star center formula at \(a\) gives \(k=m+l\);
at \(b\) it gives \(l=k+n\). Hence \(n=-m\) and \(d=-c\).
The outer polygon is nondegenerate, so \(m\ne0\).
Replace \(mT_0\) by \(T\); then \(k-l=1\), and with
\(\tau=k+l\) the centers have the form \eqref{eq:paircenters}.
The common four-star compatibility equation is
\[
 (1+12h^2)q(U)+(3+\tau^2)q(T)=0.
\]

The four outer-face determinants divided by \(D\) are
\[
 \frac{1-\tau}{4}\mathbin{\pm}\frac{h(3+\tau)}2,
 \qquad
 \frac{1+\tau}{4}\mathbin{\pm}\frac{h(3-\tau)}2;
\]
the two shared-face factors are \(1/2\pm h\tau\).
The first four are positive exactly when
\(|h|<1/6\) and \(|\tau|<(1-6|h|)/(1+2|h|)\);
these inequalities also make the last two positive.
These establish all the required necessary identities and bounds.
\end{proof}

\subsection*{Four-valent neighbors force degree six}
\begin{theorem}\label{thm:fourcore}
In every positive midpoint rectangle, each interior vertex has degree
at least four.
Every paired block has an adjacent interior vertex of degree at least
six, where a paired block consists of two adjacent interior degree-four
stars. In particular, when all interior degrees are at most five,
the degree-four vertices are pairwise nonadjacent.
\end{theorem}
\begin{proof}
At an interior vertex the nonzero, mean-zero half-eigenfunction
\(f_u=x_u-x_v\) has both signs. Its equation forces every nonzero
coordinate to have a same-sign neighbor. There are therefore at
least two coordinates of each sign, proving the degree bound.

Two adjacent interior degree-four stars necessarily have the six
faces preceding Proposition~\ref{prop:pair}. Their only possible
additional identification is \(c=d\). In that case the six faces
form a closed two-dimensional component: every edge already belongs
to two of them. Connectedness of the triangle dual of a disk excludes
such a component. Thus their boundary vertices are distinct, and
Proposition~\ref{prop:pair} makes their image a strictly convex
parallelogram \(P=(u,c,v,d)\), centered at \(o\).

If \(u\) were interior of degree four, its existing link path
\(c,a,b,d\) would have to close by the face \(udc\). But \(cd\)
is a diagonal inside \(P\), incompatible with the embedded patch.
The same holds for \(v\).

It remains to locate the common neighbors. Use the notation of
Proposition~\ref{prop:pair}, write \(s=|h|\), and express
\(T=(\kappa U_x,-\kappa U_y)\). Since \(q(U)\ne0\), this is
valid and
\[
 \kappa^2=\frac{1+12h^2}{3+\tau^2}>(s+\tfrac12)^2,
 \qquad
 3+\left(\frac{1-6s}{1+2s}\right)^2
       =\frac{4(1+12s^2)}{(1+2s)^2}.
\]
The strict inequality follows from \eqref{eq:pairdomain}.
The two coordinates of \(c-o=T-hU\) have absolute values strictly
larger than the corresponding coordinates of \(U/2\).
Thus \(u,v=o\pm U/2\) lie strictly inside the axis-parallel box
with opposite corners \(c,d\). As this box lies in the global
rectangle, both vertices are interior, hence have degree at least
five.

For the stronger degree assertion translate \(o\) to zero and put
\(k=|\kappa|\), \(r=k^2-h^2\); thus \(k>s+1/2\).
Only the case where both \(u,v\) have degree five needs proof.
The link path \(c,a,b,d\) at \(u\) is then completed by one
additional vertex \(w\), giving the two exterior faces \(uwc,udw\).
The difference of their equations is \(B(w+u,c)=0\). Substitution
in either equation gives
\[
 w=(L-\tfrac12)U-\frac{hL}{k^2}T,\qquad
 Q(L):=\frac{r}{k^2}L^2-\frac32L+\frac34-r=0.
\]
Their determinants divided by \(D\) are
\(L(r-h/2)/k^2-1\) and \(L(r+h/2)/k^2-1\).
Positivity forces \(L>k^2/(r-s/2)\ge1\).
Since \(Q(1)=1/4-k^2+h^2-h^2/k^2<0\), this selects the
unique larger root. The same calculation at \(v=-u\) therefore
forces its additional neighbor to be \(z=-w\).

Set \(A=k-s>1/2\), \(B=k+s\), and
\(L_*=k+k/(2A)\). Direct simplification gives
\[
 Q(L_*)=\frac{(A-1/2)(B-3A)}{4A}<0,
\]
because \(k>2s\). Hence \(L>L_*\). In either coordinate, put
\(\eta=\pm h/\kappa\) with the corresponding sign. The absolute
coefficient of \(c\), relative to that of \(U\), is \(k(1-\eta)\);
the coefficient of \(w\) is \(L(1-\eta)-1/2\). Since
\((L-k)(1-\eta)>[k/(2(k-s))](1-s/k)=1/2\), the latter
is positive and strictly larger. Thus \(c,d\) lie strictly inside the
axis-parallel box with opposite corners \(w,z\); they too are
global interior vertices.

The four exterior faces and the original patch fill the convex
parallelogram \((w,c,z,d)\): indeed
\(\det(w,c)/D=Lr/k^2-1/2>0\).
The vertices \(w,z\) are distinct from all old patch vertices,
since their \(U\)-coefficients are respectively greater than
\(1/2\) and less than \(-1/2\).
At \(c\) the existing link is the five-vertex path \(w,u,a,v,z\).
Degree five would close it by the edge \(wz\), a diagonal through
the already filled parallelogram, contradicting the embedding.
Thus \(c\) has degree at least six, and the same argument applies
to \(d\). This proves the stronger assertion.
\end{proof}

\section{Boundary caps and connected cores}
Delete independent interior four-stars, retaining their quadrilateral
boundaries and identities \eqref{eq:quad}. All other cells are positive
\(E=0\) triangles. Edges and vertices are inherited from the original
triangulation. In the rectangle realization these cells form an embedded
decomposition, with strictly convex quadrilaterals.

\begin{lemma}[Axial-edge rigidity]\label{lem:axialquad}
Suppose a positive four-star quadrilateral has cyclic vertices
\(A=(0,0),B=(0,L),C=(x',y'),D=(x,y)\), with \(L,x,x'>0\).
Then
\[
 x'=x,\qquad y'=L-y,\qquad -L/2<y<L/4,
\]
and its center is
\(z=(x(L-4y)/(2(L-y)),L/2)\).
Conversely these conditions give a positive four-star. Thus an axial
edge forces the opposite edge to be parallel, with the same midpoint
coordinate along that axis. A diagonal of a positive four-star cannot
be axial.
\end{lemma}
\begin{proof}
The face on \(AB\) gives \(z=(X,L/2)\), and diagonal orthogonality is
\(x'(y-L)+xy'=0\). The diagonal-intersection parameters are
\(L/(2y')\) on \(AC\) and \(L/(2(L-y))\) on \(BD\).
Strict convexity gives \(y'>L/2\) and \(y<L/2\), justifying division.
The faces on \(AD\) and \(BC\) give
\[
 X=\frac{x(L-4y)}{2(L-y)}
  =\frac{x'(4y'-3L)}{2y'}.
\]
Substituting the diagonal identity forces \(y+y'=L\) and \(x'=x\).
Strict interiority is exactly \(0<X<x\), giving the stated interval.
Conversely these formulas give all four zero face expressions and an
interior center; reflection about height \(L/2\) checks the paired faces.
If a diagonal were axial, its \(B\)-orthogonal diagonal would be
parallel, contradicting nonzero diagonal determinant.
\end{proof}

\begin{lemma}[Ordered genuine cross-sections]\label{lem:realchords}
In any positive midpoint rectangle, every original interior edge with
both endpoints on the boundary joins relative-interior points of opposite
sides. All such edges join the same pair of sides, with weakly consistent
endpoint order, and cutting them all gives a path of regions. At most
two share a boundary endpoint; if two do, their intervening region is
one original triangle with no core. No interior-degree bound is needed.
Here the cross-sections are original edges, and their slopes are unrestricted.
\end{lemma}
\begin{proof}
Suppose a genuine edge joins adjacent sides at \(U=(u,0),V=(0,v)\),
where \(u,v>0\). The original face on its corner side has third point
\((a,b)\) inside \([0,u]\times[0,v]\). But
\[
 E(U,V,(a,b))=2ab-va-ub-uv\le -uv<0,
\]
since \(va+ub-2ab=a(v-b)+b(u-a)\ge0\). This also rules out corner
endpoints, by choosing an appropriate adjacent-side description of a
proper chord. Edges joining different pairs of opposite sides would
cross internally. Thus, up to axis exchange, all are
\([(0,a_i),(W,b_i)]\), and nonintersection gives
\((a_i-a_j)(b_i-b_j)\ge0\). Distinct edges have strictly ordered
midline intersections, even when one endpoint is shared. The regions
between consecutive cross-sections therefore have a path as their dual.

For the last assertion, take a shared endpoint \(U\) on \(x=0\).
In its ordered link the neighbor abscissas satisfy
\(f_0=f_k=0\) and \(0<f_i\le W\) for \(0<i<k\). The positive
weights in \eqref{eq:M} give, on each interior radial edge,
\[
 t_{i-1}f_{i-1}+t_if_{i+1}
       =\tfrac12(t_{i-1}+t_i)f_i .
\]
No positive slot can be a weak local minimum. The sequence strictly
increases and then decreases, with at most a two-slot plateau at its
maximum. Opposite-side neighbors have value \(W\), so there are at
most two, occupying consecutive slots. They form an original face
with \(U\), whose opposite edge is physical. In particular, if
\(V_j=(W,y_j)\), then \(y_U=(y_1+y_2)/2\).
For horizontal sides, rotate the chart and reconstruct its positive
reciprocal weights; this geometric argument changes no claimed
original-mask spectral condition.
\end{proof}

\begin{lemma}[Triangular-cycle rigidity]\label{lem:trianglerigid}
Every three-cycle of original edges in a positive midpoint rectangle
is an original face. There is no interior-degree restriction.
\end{lemma}
\begin{proof}
Suppose a three-cycle \(ABC\) encloses a nonempty triangulated region.
At each of its corners the inward fan has at least one intermediate
neighbor: otherwise \(ABC\) itself would be an original face filling
the whole triangular region. Every intermediate neighbor is strictly
inside that geometric triangle, and every radial edge has both incident
faces inside it, so it inherits the same positive midpoint equation.

Translate the least corner abscissa to zero and denote the greatest
by \(b>0\). If two corner abscissas coincide, use \(x\) or \(b-x\)
at the unique extreme corner. Its inward fan has both endpoint values
\(b\), all intermediate values strictly between zero and \(b\), and
therefore a positive weak local minimum. This contradicts the
half-average recurrence in Lemma~\ref{lem:realchords}.

Otherwise order the corners so that \(x_A=0<x_B=a<x_C=b\).
At \(A\), the inward neighbor sequence has endpoints \(a,b\) and
all intermediate values in \((0,b)\). It must strictly increase:
any descending step or plateau, followed by the strictly largest
endpoint \(b\), would force an interior weak local minimum.
Let \(f_1\) be its first intermediate value. In the midpoint equation,
the following value exceeds \(f_1\), so positivity of both weights
forces \(a<f_1/2\). Hence \(b>2a\).
At \(C\), apply the same argument to \(b-x\), whose endpoint values
in increasing order are \(b-a,b\). It gives \(b>2(b-a)\), contradicting
\(b>2a\). Thus the enclosed region is a single original face.
\end{proof}

\begin{corollary}[A flag necessity]\label{cor:flag}
Every relative disk admitting positive weights satisfying the vertex
conditions \(V\) is a flag simplicial complex, with no degree restriction.
Here flag means that every clique of the original edge graph is a face.
\end{corollary}
\begin{proof}
Proposition~\ref{prop:rectangle} gives a positive midpoint rectangle,
and Lemma~\ref{lem:trianglerigid} fills every three-clique by its original
face. A plane straight-line embedding of \(K_4\) has one vertex inside
the triangle on the other three, which the same lemma forbids.
There are consequently no cliques of order at least four.
\end{proof}
Flagness supplies the first combinatorial obstruction to positive
midpoint weights.

\begin{proposition}[Two-extrema elimination]\label{prop:extrema}
Let \(S\) be a nonempty set of original interior vertices in a positive
midpoint rectangle. Some \(v\in S\) satisfies
\[
 \deg_{C[S]}(v)\le \deg_C(v)-2.
\]
If \(|S|\ge2\), there are two distinct such vertices. Consequently the
graph induced by interior vertices of original degree at most five
is \(3\)-degenerate: every nonempty induced subgraph has a vertex
of degree at most three. More generally, whenever an extremal choice
attains equality, its neighbors in \(S\) induce a path on
\(\deg_C(v)-2\) vertices.
\end{proposition}
\begin{proof}
Choose \(v\) minimizing \(x\) in \(S\). Its embedded full star surrounds
\(v\), so some neighbor difference \(f_i=x_{w_i}-x_v\) is negative.
The radial equation is
\[
 t_{i-1}f_{i-1}+t_if_{i+1}
       =\tfrac12(t_{i-1}+t_i)f_i .
\]
A negative slot must have a negative neighbor. At least two distinct
neighbors therefore have smaller abscissa and lie outside \(S\).
This proves the bound, including ties among the minimizing vertices.
If \(x\) is nonconstant on \(S\), a maximizing vertex similarly has
two strictly larger neighbors outside \(S\), giving a second choice.
Otherwise, for \(|S|\ge2\), injectivity makes \(y\) nonconstant on \(S\).
The rotated conjugate chart \((y,-x)\) has positive reciprocal weights;
its extrema give the two choices. This argument changes no original
spectral marking.

In the equality case, the two outside slots are precisely
the negative slots at a minimum, or the positive slots at a maximum.
They are adjacent in the link; otherwise both would be isolated.
The remaining slots are consecutive, so their path edges are present.
Any extra edge between them would make a three-cycle with \(v\);
Corollary~\ref{cor:flag} would make it a face and hence a link edge.
Nonconsecutive slots have no such link edge, proving the assertion.
\end{proof}
This elimination constraint identifies the core configurations to be
examined. The geometric reductions below realize the required
eliminations by positive replacement triangulations.

\begin{proposition}[Core-block interior bound]\label{prop:blockbudget}
In a positive midpoint rectangle, suppose every original interior
vertex has degree at most five.
Let \(B\) be any maximal two-connected block of \(G_5\) with at
least three vertices, and let \(Z\) be its outer facial cycle.
The original triangulated disk \(D\) bounded by \(Z\) contains
exactly the cores of \(B\) and complete four-stars. If \(I\) counts
its strictly interior cores and \(Q\) its four-star centers, then
\begin{equation}\label{eq:blockbudget}
 2Q+I+\sum_{v\in Z}(4-\deg_D(v))=6,
 \qquad 2Q+I\le2.
\end{equation}
The core graph need not be connected. This bounds interior complexity,
not the number of boundary cores or the length of a reduction.
\end{proposition}
\begin{proof}
Work first in the triangle--quadrilateral cellulation obtained by
deleting the independent four-star interiors. A simple core cycle
encloses no physical-boundary vertex. The skeleton \(H\) of the
cellulated disk inside it is two-connected: after removing a vertex,
all its neighbors remain joined by the rest of its incident cell
boundaries, in cyclic order at an interior vertex and fan order at a
boundary vertex. Every component of \(H\) minus that vertex contains
one of those neighbors, so the remainder is connected.
If the cycle lies in \(B\), then \(H\cup B\) is two-connected,
since the two graphs share a cycle. Maximality puts \(H\) in \(B\).
A two-connected plane graph has simple facial boundaries; applying
this argument to \(Z\) and to each bounded face shows that its filled
core graph is exactly \(B\), with one triangle or full four-star
quadrilateral in each bounded face.

Restore the four-star centers to obtain \(D\). Its interior degrees
are four at these centers and five at its cores. For a triangulated
disk with \(n\) interior vertices, \(b\) boundary vertices and \(e\)
edges, Euler's formula and face incidence give \(e=3n+2b-3\).
Subtracting the degree sum from \(6n+4b\) gives the identity in
\eqref{eq:blockbudget}.

Every vertex \(v\) of \(Z\) has a neighbor outside \(D\).
Indeed the exterior original triangle along a cycle edge \(vr\)
has its third vertex outside \(D\): an interior vertex cannot be
reached without crossing \(Z\), and a third vertex on \(Z\) would
force an exterior chord between vertices of \(B\). Such a chord
belongs to the maximal block, contradicting its outer boundary.
If all three edges were on \(Z\), the triangle would occupy its
bounded side, not its exterior. Hence \(\deg_D(v)\le4\), and every
boundary term \(\kappa(v)=4-\deg_D(v)\) is nonnegative.

All vertices of \(D\) are in the coordinate ranges of \(B\), since
each restored center is inside its convex quadrilateral of core
corners. The sign argument of Proposition~\ref{prop:extrema}, in
each coordinate, requires two strictly smaller and two strictly
larger original neighbors at every interior vertex of the rectangle.
Thus each coordinate extreme of \(B\) is on \(Z\), has at least
two neighbors outside \(D\), and contributes \(\kappa\ge1\).

A simultaneous extreme of both coordinates contributes at least two.
Reflect and translate so that this vertex is \(v=0\) and \(D\)
lies in the closed northeast quadrant; these changes preserve
\(E=0\). If \(\deg_D(v)=3\), its only two outside neighbors must
both be strictly southwest, by the two sign requirements.
Take the exterior original face \(vrs\) along a cycle edge \(vr\).
Here \(r\ne0\) is northeast and \(s\) is strictly southwest, so
\[
 E(0,r,s)=2q(r)+2q(s)-B(r,s)>0,
 \qquad q(w)=w_xw_y,
\]
contradicting its face equation. Consequently \(\deg_D(v)\le2\).
Both coordinate ranges are nonconstant, because \(B\) contains an
embedded cycle. Assign the four labels \(x_{\min},x_{\max},
y_{\min},y_{\max}\) to attaining vertices. Each vertex carries
at most two labels and contributes at least their number to
\(\sum\kappa\). This sum is at least four, including tied extrema,
and the identity proves the bound.
\end{proof}
\begin{corollary}[Four possible block interiors]\label{cor:blockinteriors}
With the notation of Proposition~\ref{prop:blockbudget},
\[
 (Q,I)\in\{(0,0),(0,1),(0,2),(1,0)\}.
\]
In particular, a block enclosing a four-star has exactly one such
star and no strictly interior core; a block without a four-star has
at most two strictly interior cores.
\end{corollary}
\begin{proof}
The integers $Q,I$ are nonnegative and satisfy $2Q+I\le2$.
\end{proof}
The Euler identity in~\eqref{eq:blockbudget} is classical. The face
equation at simultaneous coordinate extrema supplies the additional
inequality: four units of the total six are required at the boundary,
leaving at most two for the interior. The four resulting interior types
organize the block analysis; the later face equations determine their
positive completions. The estimate is independent of the length of the
boundary cycle.

\begin{lemma}[Corner propagation]\label{lem:corner}
A four-star quadrilateral with two consecutive sides on perpendicular
coordinate axes is an axis-aligned rectangle. If its opposite vertex
\(v\) is interior to the containing rectangle, its other two boundary
vertices are the strict side projections of \(v\). A triangle at \(v\)
cannot continue from either projection along that same boundary side.
\end{lemma}
\begin{proof}
Lemma~\ref{lem:axialquad} makes each pair of opposite sides parallel,
so the quadrilateral is an axis rectangle.
For a triangle with an axial boundary base, \(E=0\) says that the
corresponding coordinate of its apex is the midpoint of its base.
A base beginning at the apex's projection would therefore have zero length.
\end{proof}

\begin{lemma}[Empty caps]\label{lem:cap}
Let a proper straight chord of the rectangle cut off a cap containing
no retained interior vertex. Then the cap is a series of axis rectangles
coming from four-stars, and its root chord joins opposite rectangle sides.
\end{lemma}
\begin{proof}
All cell vertices in the cap lie on the physical rectangle boundary.
Its weak dual is a tree by Euler's formula. Root it at the chord.
A nonroot triangular leaf has a private vertex of original degree two,
contrary to Corollary~\ref{cor:degree}. A quadrilateral leaf has two
private degree-three vertices, both rectangle corners. A proper chord's
open boundary arc contains at most three corners, so there is at most
one nonroot leaf. Thus the dual is a path with the root at one end;
the one-cell case has the same conclusion and cannot be triangular.

The far-end quad has successive vertices
\((0,A),(0,0),(W,0),(W,B)\), up to symmetry.
Diagonal orthogonality gives \(W(B-A)=0\), making it a full-width
rectangle. Suppose the cells processed so far fill
\([0,W]\times[0,s]\). A next triangle would have perpendicular
axial legs, contradicting \(E=0\). A next quad with opposite
interfaces has vertices \((0,s),(W,s),(W,t),(0,r)\); its identity
gives \(t=r\). If the interfaces are adjacent, its two physical
edges instead run through a top corner. In the order
\((0,s),(W,s),(W,H),(a,H)\), the identity is \(a(H-s)=0\).
Then its outgoing edge lies on the physical side, contradicting its
role as an interior interface or the proper root chord. Induction
proves the assertion. No new original vertex can lie inside an edge.
\end{proof}

\begin{proposition}[The core exterior]\label{prop:core}
Assume the retained interior vertices form a nonempty connected graph
of original edges. The union \(K\) of their incident cells is a rectangle.
Its exterior consists of at most two parallel end strips as in
Lemma~\ref{lem:cap}. No degree bound is required.
\end{proposition}
\begin{proof}
An interface edge \(e\) between \(K\) and an exterior component has
both endpoints on the physical boundary: an edge incident to a retained
interior vertex has both cells in \(K\). Thus \(e\) is a proper chord.
The connected original core graph cannot cross \(e\), so all core
vertices are on one side, and every cell of \(K\) stays on that side.
The opposite cap contains no core and is exactly that exterior component;
its dual is connected and its only interface with \(K\) is \(e\).
Apply Lemma~\ref{lem:cap}. Two perpendicular end strips would overlap
at a corner, and two at the same end would overlap. Thus at most two
parallel end strips occur. Their removal leaves precisely a rectangle
filled by \(K\), proving also that its boundary is simple.
\end{proof}

We shall also use a boundary observation. Three consecutive physical
edges supporting positive \(E=0\) triangles with the same interior
apex must, up to symmetry, have vertices
\begin{equation}\label{eq:threebases}
 (0,A),(0,0),(W,0),(W,A),\qquad v=(W/2,A/2).
\end{equation}
Two consecutive bases on the same side would demand conflicting
midpoints, so the two intermediate vertices are rectangle corners;
the three midpoint equations give the formula. Four such consecutive
bases with distinct endpoints are impossible: a fourth either continues
along the same side or forces the final endpoint to equal the first.

\section{Geometric reductions and terminal rectangles}

\subsection{Local geometric reductions}\label{subsec:localreductions}

Throughout this subsection the midpoint rectangle is positive, every
interior degree is at most five, and the graph $G$ of degree-five
interior vertices is connected. We use the quadratic form $q$, its
polarization $B$, and the face defect $E$ defined above.
All triangles mentioned as faces are faces of the triangulation before
the indicated flip. A quadrilateral always denotes a complete positive
four-star, not two triangles artificially drawn across its diagonal.

\begin{figure}[htbp]
\centering
\includegraphics[width=.92\linewidth]{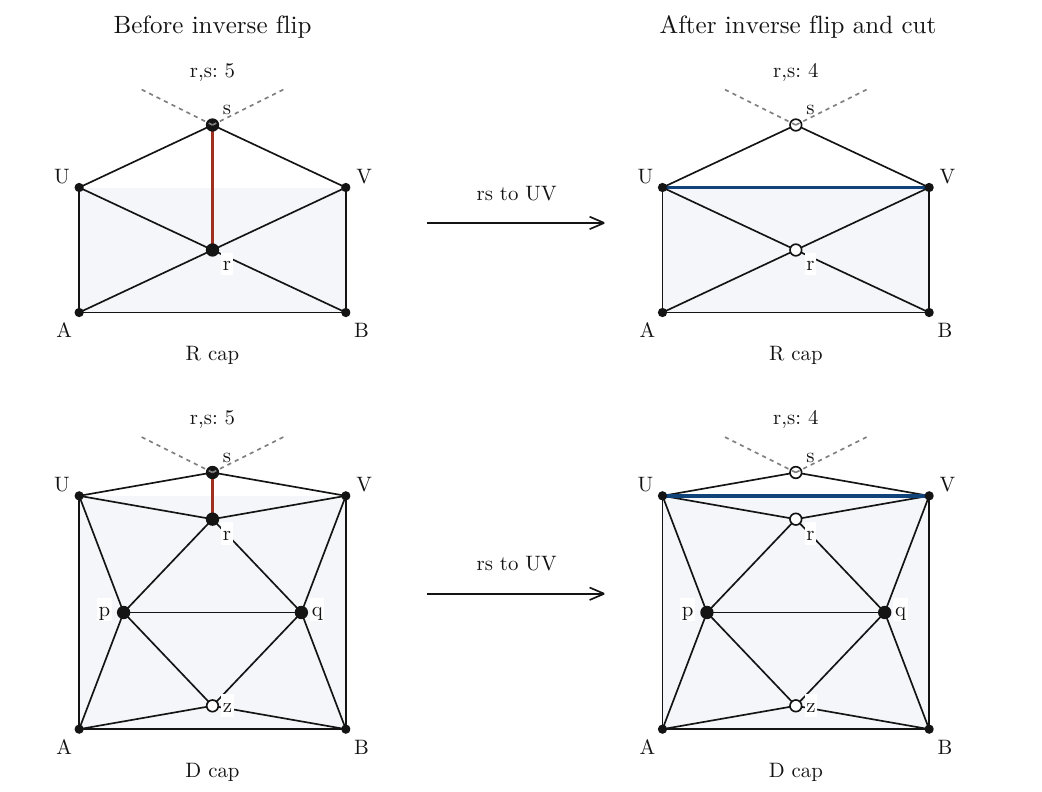}
\caption{The two inverse flips, with bridge endpoints named $r,s$.
The red edge $rs$ becomes the blue single-edge cut $UV$; both apices
drop from degree five to four. Shading locates the R or D cap after
the flip, and dashed rays indicate the omitted remainder. In the lower
row the cap also contains the two degree-five vertices $p,q$.
The original old/new labels are retained throughout the geometric cut.}
\label{fig:reductions}
\end{figure}

The coordinate eliminations below use one common linear rule. Expanding
$E(z,a,b)=2q(z)-B(z,a+b)+E(0,a,b)$ gives
\begin{equation}\label{eq:lineardifference}
\begin{aligned}
 E(z,a,b)-E(z,c,d)
 &=B(z,c+d-a-b)\\
 &\quad+E(0,a,b)-E(0,c,d).
\end{aligned}
\end{equation}
Thus differences of face equations are linear in the shared unknown
vertex $z$. We first solve these linear equations, then impose the
remaining original face equations and positive orientation. The
one-center and two-center cases use this same elimination; their
separate incidence and strict-sign checks are still necessary. The
two resulting inverse flips are shown in Figure~\ref{fig:reductions}.

\begin{lemma}[Strict coordinate slots]\label{lem:strictslots}
At every interior vertex, each coordinate has at least two strictly
larger and at least two strictly smaller neighboring values.
If $ps$ is a bridge of $G$, then $p$ has at most two other core neighbors;
if there are two, they are adjacent.
\end{lemma}
\begin{proof}
For the abscissa differences $f_i$ around the vertex, the positive radial
equations are
\[
 t_{i-1}f_{i-1}+t_i f_{i+1}
       =\tfrac12(t_{i-1}+t_i)f_i.
\]
Summing gives a positive weighted mean zero. The differences are not
all zero, since the star is embedded with positive face areas. A unique
positive entry is impossible in its displayed equation, and the same
argument applies to a unique negative entry. The conjugate positive
weights give the identical conclusion for ordinates.
For the bridge assertion, the two neighbors of $s$ in the five-cycle
link of $p$ cannot be cores: either would give a core triangle on $ps$.
The remaining two slots are consecutive, so if both contain cores they
are joined by their actual link face.
\end{proof}

We first record two exposure facts used below. A four-star with opposite
cores and two boundary corners is impossible: its boundary-to-boundary
diagonal separates the two cores, whereas a connecting core path cannot
cross the cell or pass through a boundary vertex. After this exclusion,
the noncore paths in a leaf or pendant sector have boundary-to-boundary
outer edges. Connectedness puts all cores on their inner side; the other
side is an empty cap. Proposition~\ref{prop:core} and Lemma~\ref{lem:cap}
therefore make precisely these paths physical, without changing retained
interior stars. We perform these justified removals before a reduction.

\begin{lemma}[Two exposed ears]\label{lem:twoexposedears}
Suppose $pqr,pqy,qrz,rpx$ are actual faces and the exclusive paths from
$y$ to $z$ at $q$ and from $z$ to $x$ at $r$ are physical two-edge paths.
Each path consists of two triangles or one private four-star. Both paths
are triangular, and
\[
 q=(y+z)/2,\qquad r=(z+x)/2,\qquad p=(x+y)/2.
\]
\end{lemma}
\begin{proof}
If the $q$-ear is private and the $r$-ear triangular, corner propagation
normalizes their consecutive physical corners to
\[
 y=(0,Y),\quad z=(t,0),\quad x=(W,X),\quad
 q=(t,Y),\quad r=((W+t)/2,X/2),\quad 0<t<W.
\]
Here $X,Y,W>0$. The equation $E(q,r,z)=0$ gives $X=Y$.
Then $E(p,q,y)=0$ gives $p_x=t/2$, and subtracting the equations for
$prx$ and $pqr$ gives $p_y=3Y/2$. But
$E(p,q,r)=-3WY/4\ne0$. Reflection treats the other mixed case.
If both ears are private, $q,r,z$ lie on the same line through their
common side projection, contradicting positivity of $qrz$.

Both ears are consequently triangular. Their two axial bases give the
first two midpoint formulas, and $E(q,r,z)=E(x,y,z)/4=0$.
Put $p_0=(x+y)/2$. The medial subdivision makes $p_0$ satisfy the other
three face equations. Subtracting these from the equations for $p$ gives
\[
 B(y-r,p_0-p)=B(x-q,p_0-p)=0.
\]
The two coefficient vectors have determinant equal, up to sign, to
$3\det(y-z,x-z)/4$, which is nonzero by the actual face $qrz$.
Nondegeneracy of $B$ proves $p=p_0$.
\end{proof}

\begin{proposition}[Leaf reduction]\label{prop:leafreduction}
If $|G|>2$ and $s$ is a leaf of $G$ with neighbor $p$, the two face-thirds
on $ps$ are boundary vertices. After the exposed empty caps are removed,
$ps$ is one diagonal of a strictly convex axial diamond with bisecting
diagonals. Its inverse flip cuts off a copy of $R$ along one complete
single-edge side. Both bridge endpoints lose one interior degree;
all other retained interior degrees are unchanged. The geometric
remainder is a positive midpoint rectangle with connected or empty
five-core graph.
\end{proposition}
\begin{proof}
The thirds on $ps$ cannot be cores, since $s$ is a core leaf. We must
exclude one or two four-centers; boundary exposure is as described above.

\emph{One four-center.}
Let its quad be $(s,p,t,a)$ and let the other face be $spv$.
Here $a,v$ are boundary vertices and $t$ is either boundary or a core.
The remaining five-link slots of $s$ give a physical path $a,c,v$,
consisting of two triangles or one private quad.
In the triangular case normalize
\[
 c=(0,0),\quad a=(1,0),\quad v=(0,1),\quad s=(1/2,1/2).
\]
The face $spv$ yields
\[
 p=(x,3x/D),\qquad D=4x-1,\qquad x>1/4.
\]
Writing $z$ for the four-center, the equations $E(s,z,a)=E(s,z,p)=0$
give the two candidates
\[
 z_+=(x+1/2,(2x+1)/(2D)),\qquad
 z_-=(D/(4x),D/(4(x-1))).
\]
Their difference equation is linear in $z$, so substitution in either
face equation is quadratic and gives exactly these candidates. At $x=1$
the second candidate has no finite solution. For $x<1$ it has negative
ordinate; for $x>1$, putting $C=4x^2-2x+1>0$, it has
\[
 \det(p-s,z_--s)=\frac{(4x^2-1)C}{8xD(x-1)}>0,
\]
the same side of $sp$ as its other actual third $v$. Thus $z=z_+$.
For the fourth corner, diagonal orthogonality first gives
\[
 t=s+\lambda(x-1,-3x/D).
\]
One of its remaining face equations then gives either the repeated
corner $t=s$, which is inadmissible, or
\[
 \lambda=\frac{C}{4x(x-1)},\qquad
 t=\left(x+\frac1{4x},\frac{(2x+1)^2}{4D(1-x)}\right).
\]
The value $x=1$ would make the quad diagonals parallel; $x>1$ makes
$t_y<0$. For $x<1$, the required orientation at $z$ is tested by
\[
 \det(t-z,a-z)=-\frac{(1-4x^2)C}{8xD(1-x)},
\]
and forces $1/4<x<1/2$; equality would put the center on a side.
If $t$ is boundary, the physical edge $at$ is impossible: $t_y>0$ and
$t_x>1$, while $a$ is either flat on the bottom or the bottom-right
corner of a width-one rectangle.
If $t$ is a core, the full link of $p$ is $(s,z,t,w,v)$.
Its first three listed neighbors have larger abscissa and $v$ has
smaller abscissa, so Lemma~\ref{lem:strictslots} forces $w_x<p_x$.
But $t$ already has the three smaller-abscissa neighbors $p,z,a$;
its two remaining neighbors include $w$ through the face $ptw$.
The same lemma forces $w_x>t_x>p_x$, a contradiction.

If the path $a,c,v$ is a private quad, normalize instead
\[
 c=(0,0),\quad a=(1,0),\quad v=(0,1),\quad s=(1,1).
\]
The actual face $spv$ is outside the filled square, so $p_y>1$, and
its horizontal base gives $p_x=1/2$. The axial side $as$ of the other
quad forces $t_x=p_x=1/2$ and $t_y+p_y=1$. Thus $t_y<0$.

\emph{Two four-centers.}
Write their quads as $(s,p,t,a)$ and $(s,b,u,p)$, with centers $z_1,z_2$.
The complete links force the faces $sab$ and $ptu$.
The labels $a,b$ are distinct, as are $t,u$, since equality would close
a four-cycle link at $s$ or $p$. Also $t\ne b$ and $u\ne a$, since
otherwise $pb$ or $pa$ would be a diagonal through the other quad.
If $t,u$ are both boundary, the completed links leave $G=ps$, excluded.
If $t$ is a core and $u$ boundary, the chain $a,b,u$ is physical.
For flat $b$, the axial-edge rule successively makes $ps$ and $at$
parallel to that physical line, putting $t$ on the boundary.
For a corner $b$, normalize
\[
 b=(0,0),\quad a=(A,0),\quad u=(0,U),\quad A,U>0.
\]
The face $sab$ and the quad based on $bu$ give
$s_x=p_x=A/2$ and $s_y+p_y=U$. The other quad then gives $t=(A,U)$.
The two quads with $sab,ptu$ fill the rectangle $batu$.
The four known neighbors of $t$ are $a,z_1,p,u$; its unique exterior
neighbor $w$ must give the faces $tuw,twa$. Their axial bases force
$w=(A/2,U/2)$ inside the filled rectangle, impossible for an exterior fan.
Reflection handles the other mixed case.
If $t,u$ are cores, each has four known neighbors and their common
exterior face on $tu$ identifies their unique fifth neighbor $w$.
Their fans give $tuw,taw,ubw$. It is distinct from all quad vertices:
any repetition closes a link or collapses a face, while $w=s$ gives
the leaf extra core neighbors. The two quads and
$sab,ptu,tuw,taw,ubw$ form a thirteen-face disk bounded by the original
three-cycle $a,b,w$, with six interior vertices. This contradicts
Lemma~\ref{lem:trianglerigid}.

\emph{Boundary thirds and the cut.}
Write the two thirds as $a,b$. The outer path of $s$ has three edges.
A private corner quad adjoining a triangle is excluded by
Lemma~\ref{lem:corner}; hence all three cells are triangles.
Equation~\eqref{eq:threebases} gives
\[
 a=(0,h),\quad b=(W,h),\quad s=(W/2,h/2),\quad W,h>0.
\]
The two bridge faces give $p=(W/2,3h/2)$. Thus $ps$ and $ab$ cross at
their common midpoint. Replacing $ps$ by $ab$ gives two positive
$E=0$ faces; below $ab$ is precisely $R$. The vertex $p$ remains
interior above $ab$ and has degree four. The cut is a full side with
no vertex in its relative interior. All other retained interior
stars are unchanged. The bridge rule in Lemma~\ref{lem:strictslots}
shows that deleting $s,p$ from the core graph leaves it connected or
empty: the at most two remaining neighbors of $p$ are adjacent.
\end{proof}

\begin{proposition}[Pendant-triangle reduction]\label{prop:pendantreduction}
Suppose $p,q,r$ form a core triangle, $p,q$ have no other core neighbors,
and $r$ has exactly one further core neighbor $s$, with $rs$ a bridge.
After exposed empty caps are removed, the two face-thirds on $rs$ are
boundary vertices. The inverse flip of $rs$ cuts off a copy of $D$
along one complete single-edge side. It preserves positive midpoint
geometry and the interior degree bound; the remaining core graph is
connected or empty. The four cores $p,q,r,s$ leave that graph.
\end{proposition}
\begin{proof}
The triangle $pqr$ is facial by Lemma~\ref{lem:trianglerigid}. Write the
five-link of $r$ as $(p,q,z_1,s,z_2)$. Each $z_i$ is boundary or a
four-center. If $z_1$ is a center its quad is $(r,q,a,s)$, and if $z_2$
is a center its quad is $(r,s,b,p)$. The fourth corners $a,b$ are
boundary: a core would give $q$ or $p$ another core neighbor, and
independence excludes another four-center.

\emph{Two four-center thirds.}
Their centers are distinct. If $a=b$, the known faces close the link
cycle $(r,z_1,a,z_2)$ at $s$, making its degree four. Otherwise $s$ has
exactly the neighbors $r,z_1,a,b,z_2$, so the face $sab$ is forced.
If the outside cell at $pq$ is a triangle $pqw$, the links force $pbw$
and $qwa$. The two quads and $pqr,pqw,pbw,qwa,sab$ form a thirteen-face
disk bounded by $awb$ with six interior vertices. All boundary labels
are distinct, since any equality repeats a link vertex or collapses a
forced face. This contradicts faciality of three-cycles.
If the outside cell at $pq$ is a quad, its corners must be $(p,q,a,b)$
by the last available link slots. Its center $z_0$ and the other faces
close the cycle $(q,z_1,s,b,z_0)$ in the link of the boundary vertex $a$.
A boundary link is a path, so this too is impossible.

\emph{One four-center third.}
By reflection let the quad be $Q=(r,q,a,s)$ and the other third be a
boundary vertex $v$. The remaining actual faces include $pqr,rpv,rsv$.
We have $a\ne v$, since otherwise $ra$ is a diagonal through $Q$.
The outside $pq$ cell is either an ordinary quad $(b,a,q,p)$ or a
triangle $pqw$. These are exhaustive: an additional core there gives
$p$ or $q$ another core neighbor, and a quad containing $p,q,r$ has
an existing core edge as its diagonal.

Here is a common elimination used to check every case without discarding
an unused equation. For $s=(X,Y)$ set
\[
 F=q(r-a)+q(q-s)+3q(r+a-q-s).
\]
The full four-star identities give
\begin{equation}\label{eq:pendantquad}
 B(r-a,q-s)=0,\qquad F=0.
\end{equation}
The combination $2E(r,s,v)-F=0$ is the affine equation
\begin{equation}\label{eq:pendantlinear}
\begin{split}
 &B(r+3a-2q-2v,s)+4q(r)+4q(v)-2B(r,v)\\
 &\hspace{12mm}-q(r-a)-q(q)-3q(r+a-q)=0.
\end{split}
\end{equation}
Thus the first equation of \eqref{eq:pendantquad} and
\eqref{eq:pendantlinear} are linear in $X,Y$. We use these necessary
equations below, then check an unused actual face or the axial-edge
four-star condition.

If the outside cell is $(b,a,q,p)$, its physical boundary chain is
$v,b,a$. The axial-edge lemma normalizes it to
\[
 b=(0,0),\quad a=(1,0),\quad p=(t,1),\quad q=(1-t,1),
 \qquad -1/2<t<1/4,
\]
and $pqr$ gives $r=(1/2,R)$. For corner $b$, the physical face $pvb$
gives $v=(0,2)$ and $0<t<1/4$. The equations above then give
\[
 R=\frac3{2(1-t)},\quad
 X=\frac1{4t},\quad Y=\frac3{4t(1-t)}-2,
\]
but the unused face satisfies
\[
 8t^2(1-t)E(r,s,v)=3(1-2t)(1-4t)>0.
\]
For flat $b$, physical order gives $v=(2t,0)$ and $t<0$. Now
\[
 R=\frac1{2(1-3t)},\quad
 X=1+2t-\frac1{4t},\quad
 Y=\frac{1-1/(4t)}{1-3t},
\]
and
\[
 2(1-3t)E(r,s,v)=\frac{(4t-1)(2t-1)}{4t^2}>0.
\]
The value $t=0$ collapses the physical flat edge or puts $p$ on a
physical side at a corner. Hence this outside-quad case is impossible.
No claim that $a_x=1$ is the global width is used.

Suppose instead the outside cell is $pqw$. The exclusive $q$-cell is
$qwa$, and the physical $p$-path is $v,c,w$. If this path is triangular,
normalize
\[
 c=(0,0),\quad v=(0,1),\quad w=(1,0),\quad p=(1/2,1/2).
\]
Writing $q=p+(x,y)$ and $r=p+(u,\zeta)$, the three actual face equations
$pqw,rpv,pqr$ give
\[
 (4x-1)y=1-x,\qquad (4u+1)\zeta=u+1,
\]
\[
 2(2x-1)u^2-(4x^2-6x+1)u-x(2x-1)=0.
\]
Its discriminant is $(4x^2-2x+1)^2$, and the two candidates are
\[
 r_1=\left(x,\frac{3x}{4x-1}\right),\qquad
 r_2=\left(-\frac1{2(2x-1)},\frac3{2(2x+1)}\right).
\]
If $w$ is flat, write $a=(A,0)$ with $A>1$ and $D=2A-1$.
The physical face $qwa$ gives
\[
 q=\left(\frac{A+1}{2},\frac{A+1}{2D}\right),\quad
 x=A/2>1/2.
\]
The second root has negative abscissa, so $r=(A/2,3A/(2D))$.
The two linear equations reduce to
\[
 Y=(3X-A-1)/D,\qquad 4D(X-A)=0,
\]
hence $s=(A,1)$. The actual side $as$ is vertical, so the opposite
side $rq$ of its positive four-star must be vertical, whereas
$q_x-r_x=1/2$.
If $w$ is a corner, write $a=(1,A)$, $A>0$, and $D=2A-1$. The actual global width
is now one, and
\[
 q=\left(\frac{3A}{2(2A-1)},\frac A2\right).
\]
Interiority gives $A>2$; $A=1/2$ is inconsistent with the face equation.
Now $x=(A+1)/(2D)\in(1/4,1/2)$ and $(r_2)_x=D/[2(A-2)]>1$,
so $r=((A+1)/(2D),(A+1)/2)$. The linear equations are
\[
 Y=A-DX/3,\qquad (7A-5)(Y-A)+D(5A-3)X=0.
\]
They give $4D^2X/3=0$, putting $s=(0,A)$ on the boundary. In both
ranges $x\ne1/2$, so no vanishing quadratic coefficient was omitted.

Finally, if $v,c,w$ is a private quad, corner propagation gives
\[
 c=(0,0),\quad v=(0,1),\quad w=(1,0),\quad p=(1,1).
\]
The vertices $v,w$ are flat because $p$ is interior. Thus $a=(A,0)$,
$A>1$, and the actual faces give
\[
 q=((A+1)/2,1/2),\quad r=(1/2,3/[2(A+1)]),
\]
\[
 X=\frac{A(A-2)}{2(A-1)},\qquad
 Y=\frac{A^2-4}{2(A^2-1)}.
\]
For $1<A\le2$, $X\le0$. For $A>2$ the unused face gives
\[
 E(r,s,v)=-\frac{3A^2(A-2)}{2(A-1)^2(A+1)}<0.
\]
At $A=2$ the putative $s$ is the physical corner, not a lost solution.
This exhausts the one-center configurations and proves that the bridge
thirds are boundary vertices.

\emph{A clean bridge.}
Write its thirds as $u,v$, with link $(p,q,u,s,v)$ at $r$.
If the outside $pq$ cell is a triangle $pqw$, the exposed ears at $p,q$
and Lemma~\ref{lem:twoexposedears} give
\[
 p=(v+w)/2,\quad q=(w+u)/2,\quad r=(u+v)/2,
 \qquad E(u,v,w)=0.
\]
Put $U_0=u-v$, $S_0=s-r$. The two actual bridge equations give
\[
 B(U_0,S_0)=0,\qquad q(S_0)=-q(U_0)/4.
\]
An axial $U_0$ makes the bridge faces degenerate. Otherwise
$S_0=\pm((U_0)_x/2,-(U_0)_y/2)$, so $s$ is one of the other two
corners of the axis box on diagonal $uv$. If $u,v$ are on opposite
physical sides both choices are boundary; the same-side case is axial.
For adjacent sides, the only potentially interior choice is, after
normalization,
\[
 u=(a,0),\quad v=(0,b),\quad s=(a,b),\quad
 0<a<W,\quad0<b<H.
\]
The two physical ears have two distinct private corners. Their four-edge
path through $w$ must therefore take the long boundary arc, giving
$w=(W,H)$ and
\[
 2WH=aH+bW+ab.
\]
For $A=a/W$, $B=b/H$, this says $2=A+B+AB$, whence $A,B>1/2$ and
$|A-B|<AB$. Consequently
\[
 \det(p-r,s-r)=\frac{Wb-ab-aH}{4}<0,\qquad
 \det(s-r,q-r)=\frac{aH-ab-Wb}{4}<0.
\]
The ray $rs$ lies strictly inside the angle of the actual face $pqr$,
contrary to its exterior bridge sector.

The outside cell is therefore the ordinary quad $(A,B,q,p)$.
Normalize its physical base to $A=(0,0)$, $B=(W,0)$. The axial-edge
four-star lemma gives
\[
 p=(t,H),\quad q=(W-t,H),\quad -W/2<t<W/4,
 \quad z=\left(W/2,\frac{H(W-4t)}{2(W-t)}\right).
\]
The central triangle gives $r=(W/2,R)$; its physical outer chain is
$u,A,B,v$. If both $A,B$ are corners,
\[
 u=(0,2H),\quad v=(W,2H),\quad0<t<W/4.
\]
If both are flat, $u=(2t,0)$, $v=(W-2t,0)$, and $t<0$.
A mixed case is impossible: the corner end requires $t>0$, whereas
the flat end requires $t<0$. In the flat case
\[
 R=\frac{WH}{2(W-3t)}>0,\qquad s_y=-R<0,
\]
by $rpu$ and the bridge equations, so only the corner case remains.
It gives
\[
 R=\frac{3WH}{2(W-t)},\quad \frac{3H}{2}<R<2H,
 \qquad s=(W/2,4H-R).
\]
Thus $rs$ and $uv$ are perpendicular axial diagonals with common
midpoint $(W/2,2H)$, interior to both. The inverse flip replaces the
two bridge faces by positive $E=0$ faces. Below $uv$ is precisely $D$:
its boundary is $(A,B,v,u)$, its old cores are $p,q$, and its two
four-centers are $z,r$. Its ten faces are
\[
 ABz,\ Bvq,\ vur,\ uAp,\ Bzq,\ vqr,\ urp,\ Apz,\ zqp,\ pqr.
\]
The vertex $s$ is above the cut and becomes degree four. No other
retained interior degree changes. The outside endpoint $s$ has at
most two other core neighbors, adjacent when both exist; deleting
$p,q,r,s$ therefore leaves the core graph connected or empty.
\end{proof}

Both reductions preserve the positive $E=0$ geometry. The isolated
remainder admits the compatible marking of its own rectangle; this
auxiliary choice does not change any original-old or original-new
label in the spectral problem. In particular no preservation of the
original empty-face positive-semidefinite condition is inferred from
these geometric propositions.

\subsection{The terminal rectangles}\label{subsec:terminalrectangles}

\begin{proposition}[Terminal bases]\label{prop:terminalbases}
For a positive midpoint rectangle of interior degree at most five whose
core graph is empty or connected and has at most three vertices:
\begin{enumerate}
\item an empty core graph gives a nonempty serial strip of $R$ pieces;
\item a nonempty triangle-free core graph of order at most three is $K_2$;
      after empty end strips are removed its support is $F$ or $D$;
\item a three-vertex core graph is $K_3$; its incident-core rectangle is
      the ten-face piece $T$; the whole rectangle may additionally have an
      unflipped $R$ strip on at most one of its two eligible adjacent sides.
\end{enumerate}
Here $F$ consists of two $R$ pieces with their common centered interface
flipped. Its only eligible attachment sides are its two unsubdivided
opposite ends. All four sides of $D$ are single edges, so its empty
strips may use either opposite pair, but cannot turn.
\end{proposition}
\begin{proof}
\emph{Zero and one core; the three-vertex path.}
Collapse the independent four-stars. If there are no cores and only
one cell, it must be a quadrilateral spanning the rectangle, hence $R$.
If there is more than one cell, an interior cell edge is an original
boundary-to-boundary chord. Both complementary caps contain no core;
Lemma~\ref{lem:cap} makes both serial rectangular $R$ strips, with
common cross-section that chord.

A sole core is impossible. After empty caps are removed,
Proposition~\ref{prop:core} makes its entire reduced fan a rectangle.
If the fan has no quad, at most one triangular base can lie on each
side: two disjoint bases cannot both have the core projection as their
midpoint. Its degree is then four, not five. If it has a quad,
Lemma~\ref{lem:corner} makes that cell the axis rectangle between the
core and a physical corner. At either strict side projection the next
cell cannot be triangular. Repeating around the four sides forces four
corner quads and therefore original degree eight, again impossible.
If $G=P_3$, Proposition~\ref{prop:leafreduction} applies to either leaf
and leaves a positive rectangle with one core, just excluded.

\emph{Two cores.}
Write the core edge as $pq$ and remove empty exterior strips.
A quad with opposite cores and two boundary vertices is excluded by
connectedness as above. Let $j\in\{0,1,2\}$ count the shared quads
on $pq$, and let $k_p,k_q$ count all quads incident to the two cores.
Independence of four-centers gives $k_p,k_q\le2$. The number of quads
is $k_p+k_q-j$ and that of triangles is $8-2(k_p+k_q)+j$.
The disk incidence identity consequently gives perimeter length $6-j$.
All perimeter vertices are distinct because this region is a rectangle.

For $j=0$, let $u,v$ be the two boundary thirds on $pq$.
Each exclusive outer path has three edges, consisting either of three
triangles or a private quad next to a triangle. The latter is impossible
by corner propagation. The four private vertices of these paths are
therefore the four rectangle corners. Normalize to $[0,W]\times[0,H]$
with $u=(0,a)$, $v=(W,b)$, $p$ above $q$. The physical bases give
\[
 p_x=q_x=W/2,\quad p_y=(H+a)/2=(H+b)/2,\quad q_y=a/2=b/2.
\]
Thus $a=b$, and the vertical-base equation on $pqu$ gives
$a=(H+2a)/4=H/2$. Hence
\[
 p=(W/2,3H/4),\quad q=(W/2,H/4),\quad
 u=(0,H/2),\quad v=(W,H/2).
\]
This is precisely $F$. The two sides containing $u,v$ cannot receive a
full single-edge strip, leaving only the opposite unsubdivided ends.

For $j=1$, let the shared triangle be $pqw$, the quad $(p,q,u,v)$,
and the physical five-cycle be $(w,A,v,u,B)$, with exclusive paths
$w,A,v$ and $u,B,w$. Their private middle points $A,B$ are corners,
whether the ears are triangular or private quads. Thus exactly one
of $w,u,v$ is flat.
If $w$ is flat, normalize the perimeter to
\[
 v=(0,0),\ u=(W,0),\ B=(W,H),\ w=(a,H),\ A=(0,H),\quad0<a<W.
\]
A private ear would put its core on the opposite physical side by
corner propagation; hence both ears are triangular. Their bases give
$p=(a/2,H/2)$, $q=((W+a)/2,H/2)$, and $pqw$ gives $a=W/2$.
The quad center has abscissa $W/2$ from its bottom base, but its side
face has $E(z,u,q)=-3Wz_y/4$, forcing $z_y=0$.
If $v$ is flat, normalize instead
\[
 A=(0,0),\ v=(a,0),\ u=(W,0),\ B=(W,H),\ w=(0,H),\quad0<a<W.
\]
Again a private ear puts its core on a physical side. Triangular ears
give $p=(a/2,H/2)$, $q=(W/2,H/2)$, whereas $pqw$ requires
$0=w_x=(a+W)/4>0$. Reflection covers flat $u$. Thus $j=1$ is impossible.

For $j=2$, both cores already have two four-center neighbors, so there
are no private quads and each core has one further triangular cell.
The perimeter has four vertices; the two exclusive triangle bases are
opposite sides, and the shared quads occupy the other pair. Reinserting
their centers gives exactly $D$ with boundary $(A,B,V,U)$ and faces
\[
 ABz,\ BVq,\ VUr,\ UAp,\ Bzq,\ Vqr,\ Urp,\ Apz,\ zqp,\ pqr.
\]
Its complete positive coordinate family, sufficient also to verify
realizability, is
\begin{equation}\label{eq:Dcoordinates}
\begin{gathered}
 A=(0,0),\ B=(W,0),\ V=(W,H),\ U=(0,H),\\
 p=(kW,H/2),\quad q=((1-k)W,H/2),\\
 z=(W/2,hH),\quad r=(W/2,(1-h)H),\\
 k=\frac{1-4h}{4(1-h)},\qquad 0<h<1/4.
\end{gathered}
\end{equation}
Indeed the two vertical triangle bases give the core ordinates, and
the axial-edge rule on the lower quad gives $p_x+q_x=W$ and
$h=(1-4k)/(4(1-k))$; the upper quad gives its reflected center.
Interiority gives $0<k<1/4$, equivalently the displayed interval for
$h$. Conversely the axial-edge four-star lemma verifies both positive
quads, and the remaining two axial-base triangles are positive with
$E=0$. Any rectangle aspect ratio is allowed. All four perimeter sides
are single edges; disjoint nonempty strips on perpendicular sides would
overlap, so exterior strips use one opposite pair.

\emph{Three cores forming a triangle.}
The core triangle $pqr$ is an actual face. No quad can contain all three
cores, because one of its diagonals would already be a core edge.
Let $j$ be the number of shared quads on its three edges. A core with
$0,1,2$ incident shared quads has an exclusive perimeter path of length
$2,1,0$, respectively. A private quad can only replace a length-two
path. Counting cell incidences at the three degree-five vertices and
using Euler's formula gives perimeter length $6-j$. Thus $j\le2$.
A private quad requires two distinct strict flat projection vertices,
so it is impossible when $j\ge1$, since the perimeter has at most five
vertices. When $j=0$, Lemma~\ref{lem:twoexposedears} excludes any private
ear using the four actual central and shared triangles.

For $j=1$ put the shared quad on $pq$. Its boundary edge is $ab$;
the exclusive edges at $p,q$ are $wa,bz$, the exclusive path at $r$
is $z,c,w$, and the other shared triangles have thirds $w,z$.
The physical boundary is the pentagon $(w,a,b,z,c)$, with $c$ a corner.
Reflection reduces the unique flat point to $w$ or $a$.
For flat $w$ normalize
\[
 c=(0,0),\ w=(t,0),\ a=(1,0),\ b=(1,1),\ z=(0,1),\quad0<t<1.
\]
The boundary equations give
$p=((1+t)/2,P)$, $q=(1/2,Q)$, $r=(t/2,1/2)$.
The two shared triangle equations give
\[
 P=\frac{1+t}{2(2-t)},\qquad Q=\frac3{2(2-t)},\qquad Q-P=1/2.
\]
The quad center $k_0$ has ordinate $1/2$ from its vertical base $ab$;
then $E(a,p,k_0)=0$ gives $(k_0)_x=1+t/2>1$.
For flat $a$ normalize
\[
 c=(0,0),\ w=(1,0),\ a=(1,t),\ b=(1,1),\ z=(0,1).
\]
Now $r=(1/2,1/2)$ and $q=(1/2,Q)$, while $E(q,r,z)=0$ requires
$Q=3/2$. Both cases contradict interiority.

For $j=2$ put the shared quads on $pq$ and $qr$, and normalize the
remaining shared triangle third and other perimeter vertices to
\[
 w=(0,0),\ a=(1,0),\ c=(1,1),\ b=(0,1).
\]
The exclusive bases give
$p=(1/2,P)$, $r=(R,1/2)$, $q=(X,Y)$ with $0<P,R<1$.
The remaining shared triangle and the two quad diagonal equations yield
\[
 P+R-PR=1/4,\qquad
 Y=2(1-P)(1-X),\qquad X=2(1-R)(1-Y).
\]
Since $(1-P)(1-R)=3/4$, these imply $X=1/2+R$, $Y=1/2+P$.
The center $k_0$ of $(p,q,c,a)$ has ordinate $1/2$ from its vertical
base. Its face $E(a,p,k_0)=0$ gives $(k_0)_x=1-R$, but
\[
 E(p,q,k_0)=R(1+2P)>0.
\]
Thus $j=2$ is impossible as well.

Only $j=0$ remains. All cells are triangles. The three private middle
vertices are three distinct corners; one shared terminal is the fourth
corner and the other two are strict flat points on adjacent sides.
Relabel them as the boundary cycle $(A,U,B,V,C,D)$, with
\begin{equation}\label{eq:Tcoordinates}
\begin{gathered}
 A=(0,0),\ U=(a,0),\ B=(1,0),\ V=(1,b),\ C=(1,1),\ D=(0,1),\\
 p=(a/2,1/2),\quad q=((1+a)/2,b/2),\quad r=(1/2,(1+b)/2),\\
 b=\frac{1+a}{2-a},\qquad 0<a<1/2.
\end{gathered}
\end{equation}
The six exterior midpoint equations give the displayed core positions.
Each of the remaining four face equations is
$(-1-a+2b-ab)/4=0$, which gives exactly the stated range.
The ten faces are
\[
 DAp,\ AUp,\ UBq,\ BVq,\ VCr,\ CDr,\ pUq,\ qVr,\ rDp,\ pqr.
\]
Their six exterior doubled areas are
$a/2,a/2,(1-a)b/2,(1-a)b/2,(1-b)/2,(1-b)/2$; all four inner doubled
areas equal $(1-a+ab)/4$. These are strictly positive, proving the
converse realization. The sides $AB,BC$ are subdivided; only $DA,CD$
are single edges. They are adjacent, so at most one can receive an
empty strip. This proves every terminal assertion.
\end{proof}

\section{Terminal core blocks}
Throughout this section the original interior degrees are at most five,
the original five-core graph $G$ is connected, and the empty exterior
caps of Proposition~\ref{prop:core} have been removed. Thus an original
boundary-to-boundary edge incident to the core region is physical.
All stars and triangular faces below are original. In particular,
the equation $E=0$ is never imposed on an artificial half-quadrilateral.
We use flagness, independence of four-centers, the strictly convex
four-star rule, and Proposition~\ref{prop:blockbudget}.

We repeatedly use the \emph{two-sign rule} of
Lemma~\ref{lem:strictslots}: an interior vertex has at least two
strictly smaller and two strictly larger neighbors in each coordinate.

\begin{lemma}\label{lem:cutslots}
If a core $a$ is a cut vertex of a nontrivial maximal two-connected
block $K$, it has exactly two neighbors $u,v$ in $K$, joined by an
edge, and exactly one other core neighbor $b$. The edge $ab$ is a
bridge. In particular a proper terminal block has one cut, which is
a degree-two ear of its filled disk.
\end{lemma}
\begin{proof}
The $K$-neighbors belong to one component of $G-a$, and a cut has
a neighbor in another component. Consecutive core neighbors in the
original five-link are joined by a facial edge. Two components thus
require two intervening noncore slots. Two-connectedness requires
two $K$-neighbors, so the link is $u,v,z_1,b,z_2$, up to orientation.
It supplies the face $auv$; $b$ is the only neighbor in its component
of $G-a$, proving bridgehood. In the filled block disk the two edges
at $a$ bound that facial ear.
\end{proof}

\subsection{Three reusable face obstructions}
\begin{lemma}\label{lem:quadrantdichotomy}
Let $a,b>0$, $t=(a,b)$, and $z=(X,Y)$. The equation $E(0,t,z)=0$ gives
\[
 (2X-a)Y=b(X-2a),\qquad
 \det(t,z)=-\frac{2b(X^2-aX+a^2)}{2X-a}.
\]
If $X,Y>0$, positive determinant is equivalent to $X<a/2$, $Y>2b$;
negative determinant is equivalent to $X>2a$, $Y<b/2$.
For the positive-determinant implication it suffices to assume $X>0$.
\end{lemma}
\begin{proof}
Expand $E$. The denominator cannot vanish, since the first equation
would then give $0=-3ab/2$. Its numerator quadratic is positive.
The signs give the assertions, including the last one. Reflections
and positive independent axis scalings preserve the assertions with
the corresponding determinant signs. Across an actual edge the two
actual face-thirds have opposite nonzero determinant signs.
\end{proof}

\begin{lemma}\label{lem:cornerobstruction}
Put $m=(0,0)$, $A=(1,0)$, $B=(0,1)$. The eight actual faces
\[
 md_2A,\ md_1B,\ md_1d_2,\ pd_1d_2,
 \ Ad_2a,\ pd_2a,\ Bd_1b,\ pd_1b
\]
cannot all have zero $E$ if $d_1,d_2,p,a,b$ lie in the open unit square.
\end{lemma}
\begin{proof}
The first three equations give
\[
 d_1=(x,\tfrac12),\quad d_2=(\tfrac12,y),\quad
 x+y-xy=\tfrac14,\quad 0<x,y<\tfrac14.
\]
For $p=(P,T)$, subtract the equations with thirds $m,p$ to obtain
\[
 T=\frac{P(1/2+y)}{2P-x-1/2}.
\]
The denominator is positive. This branch is strictly decreasing
and passes through $(P_0,T_0)=(x+1/2,y+1/2)$.
Writing $a=(X,Y)$, its boundary equation is
$(4X-3)Y=2Xy$, so $X>3/4$ and $Y>2y$.
The algebraic comparison point $a_0=2d_2=(1,2y)$ satisfies
$E(A,d_2,a_0)=0$. Subtraction gives
\[
 E(p,d_2,a)-E(p,d_2,a_0)
       =(1-P)(Y-2y)-T(X-1)>0.
\]
Thus $E(p,d_2,2d_2)<0$. Expanding it yields
\[
 (\tfrac12-2y)(P-P_0)+(x-1)(T-T_0)<0,
\]
which forces $P<P_0$, $T>T_0$. Interchanging the coordinates, the
two actual faces at $b$ similarly imply
\[
 (y-1)(P-P_0)+(\tfrac12-2x)(T-T_0)<0.
\]
Both terms are now positive, a contradiction. The comparison points
are used only algebraically, not as additional faces.
\end{proof}

\begin{lemma}\label{lem:twobridges}
Let $W=(A,p,B,q)$ be a simple original four-cycle, with $A,B$ physical
boundary vertices and $p,q$ cores. Suppose its filled original disk
$D$ has $\deg_D(p)=\deg_D(q)=4$. Then the fifth neighbors $h,f$ of
$p,q$ are cores, and $ph,qf$ are distinct bridges with distinct
outside endpoints.
\end{lemma}
\begin{proof}
The disk has an interior vertex. A diagonal of $W$ would, by flagness,
make its two halves actual faces and leave no interior vertex.
Consequently the full links have orders
\[
 A,u,v,B,h\quad\hbox{and}\quad B,s,t,A,f,
\]
where $u,v,s,t$ are strictly inside $W$, while $h,f$ are outside.
If $h$ is boundary, the physical edges $Ah,hB$ cannot be collinear:
their midpoint equations would force $A=B$. Normalize their true
corner as $h=0,A=(1,0),B=(0,1)$, so $p=(1/2,1/2)$.
The two-sign rule forces $q\in(0,1)^2$: if, for example, $q_x\ge1$,
its two strictly interior disk neighbors are strictly left and $A,B$
are weakly left, leaving only $f$ for two right slots. Every inner
point is therefore in the unit square. But $v=(X,Y)$ satisfies
$(4X-1)Y=3X$, which for $0<X<1$ and $Y>0$ gives $Y>1$.

If $h$ is a four-center, its quad is $(p,A,l,B)$. Independence and
the connected-core obstruction to opposite cores across a BB diagonal
make $l$ boundary. Its two physical sides meet at a true corner;
the four-star rule gives
\[
 l=0,\quad A=(a,0),\quad B=(0,b),\quad
 p=(a,b),\quad h=(a/2,b/2),\quad a,b>0.
\]
At $p$ both inner neighbors must have abscissa $>a$. Hence $q_x>a$,
since those points lie inside $W$. Now $q$ is the unique outer
abscissa maximum, and again only $f$ can be to its right. Thus $h$
is a core; the reflected proof treats $f$.
The original crosscut $A-p-B$ separates $h$ from the inner disk.
A core path cannot cross an original edge or pass through $A,B$;
its only crossing is $p$, and $ph$ is its only edge into that exterior
side. Hence $ph$ is a bridge. The other crosscut gives $qf$; their
outside regions are disjoint, so $h\ne f$.
\end{proof}

\subsection{The exact block templates}
For a nontrivial maximal block $K$, write $D$ for its filled disk,
$C$ for its outer cycle, $Q$ for its enclosed full four-stars and
$I$ for its strictly interior cores. Proposition~\ref{prop:blockbudget}
gives $2Q+I\le2$. Its proof also gives $\deg_D(v)\le4$ on $C$ and
degree five at every strictly interior core.

\begin{lemma}\label{lem:noquadterminal}
A block containing a full all-core four-star has two bridge exits.
It is neither proper terminal nor the entire two-connected $G$.
\end{lemma}
\begin{proof}
Let its cyclic corners be $c_0,\ldots,c_3$, with center $z$.
Each corner has two exterior slots. If $t_i$ is the exterior third
on $c_ic_{i+1}$, the full stars force the eight faces
$c_ic_{i+1}t_i,c_it_{i-1}t_i$. Adjacent $t_i$ are distinct by the
slots; $t_i=t_{i+2}$ would give a third face on $c_{i+1}c_{i+2}$.
They differ from the inner vertices. The twelve actual faces form
a disk with simple outer cycle $(t_0,t_1,t_2,t_3)$ and no diagonal.
Each $t_i$ has four known neighbors; it cannot be a four-center,
whose closing link would supply a forbidden diagonal.

Adjacent outer cores, say $t_0,t_1$, have a common fifth neighbor
$h$ and force $ht_3t_0,ht_0t_1,ht_1t_2$. Here $h$ is new: the
old inner links are complete and each other coincidence repeats a
slot or makes a diagonal. These faces leave the separating original
triangle $ht_2t_3$, contrary to flagness. Hence the outer cores form
an independent set of the four-cycle.
If there are none, its four physical sides form a rectangle. The
four corner-base midpoint equations make the two inner-quad diagonals
horizontal and vertical, contradicting their $B$-orthogonality.

If there is one, relabel the inner quad $(v,p,u,q)$ and outer cycle
$(A,B,n,m)$, with $A,B,m$ boundary. A flat $A$ would leave $n$ as
the unique outer point off one side and only its fifth neighbor
available in the outward direction. Thus normalize
$A=0,B=(0,1),m=(1,0)$. The same sign argument bounds $n$ and all
inner vertices in $(0,1)^2$. The actual faces give
\[
 v=(x,\tfrac12),\ q=(\tfrac12,k),\quad
 k=\frac{1/4-x}{1-x},\quad 0<x,k<\tfrac14.
\]
For $p=(P,T),u=(U,V)$ their equations are
$(2P-x)T=3P/2$ and $(4U-3)V=2Uk$, so $T,U>3/4$.
Consequently
\[
 B(v-u,p-q)=(x-U)(T-k)+(\tfrac12-V)(P-\tfrac12)
                    <-\tfrac14+\tfrac14=0,
\]
contrary to the full four-star identity. The only remaining case
has two opposite outer cores. Lemma~\ref{lem:twobridges} applies.
The inner four cores and these two ear tips form a two-connected
six-vertex graph, every inner star is complete, and its only exits
are those two bridges. It is therefore the entire maximal block.
\end{proof}

\begin{lemma}\label{lem:blocktemplates}
For a proper terminal block, or $G$ itself when two-connected, its
filled disk is purely triangular and has one of the following forms:
\begin{enumerate}
\item $I=0$: the strip $P_n^2$, with faces $(i,i+1,i+2)$, $n\ge3$,
or a triangle with one single ear on each side;
\item $I=1$: a five-wheel with zero, one or two ears on a matching
of its rim;
\item $I=2$: the union of two adjacent full five-stars, with at
most one ear on each of two specified disjoint boundary edges.
\end{enumerate}
\end{lemma}
\begin{proof}
The preceding lemma excludes $Q>0$, so $I\le2$.
We first justify the attachment rule used below. Let a triangulated
subdisk $S$ have its boundary vertices on $C$, each of degree three
or four in $S$. An edge of $\partial S$ not on $C$ is a crosscut.
Its outside patch is bounded by that edge and the intervening $C$-arc,
which contains no other vertex of $\partial S$. These patches have
disjoint interiors. The first outside triangle on $uv$ introduces
a new neighbor at each endpoint; neither can already have degree
four in $S$. Both are then saturated at four. A further triangle
across either new edge would add another neighbor at its saturated
endpoint. Thus the patch is exactly one ear, and occupied edges form
a matching. In particular no strict interior vertex can be added.
The same argument gives no extension if every boundary edge touches
a saturated vertex, irrespective of its other endpoint's degree.

For $I=1$, the full star of the interior core is a five-wheel whose
rim lies on $C$, with rim degrees three. Apply the attachment rule.
For $I=2$, nonadjacent interior cores are impossible: the full star
of either would have all its rim on $C$ and leave no place for the
other. For adjacent $o,p$, flagness gives exactly two common neighbors
$a,b$. Write their complete links
\[
 a,c_1,c_2,b,p\quad\hbox{and}\quad b,d_1,d_2,a,o.
\]
The six noncentral vertices are distinct: a repeated $c$- and
$d$-vertex would be a third common neighbor. The star union has
boundary $(a,c_1,c_2,b,d_1,d_2)$, degree four at $a,b$ and degree
three at the other four. The only allowed ear edges are $c_1c_2$
and $d_1d_2$.

For $I=0$, the weak face dual is a tree of maximum degree three.
A branching face $abc$ has exterior thirds $x,y,z$ on its three
sides. They are distinct: two equal thirds would close a link and
make a boundary vertex interior. The four triangles form a six-vertex
disk; $a,b,c$ are saturated, and every boundary edge meets one of
them. The attachment rule permits no extension.
Otherwise the dual is a path. Removing its end triangle removes a
private degree-two boundary corner, so in reverse order each triangle
introduces a new vertex. Label the first three faces $123,234,345$,
interchanging $2,3$ if necessary. Vertex $3$ is now saturated, so
the next face must be $456$, not a face on $35$. Inductively the
saturated older corner forces the consecutive-triple strip. The
one- and two-face cases are $P_3^2,P_4^2$.
\end{proof}

\begin{figure}[htbp]
\centering
\includegraphics[width=\linewidth]{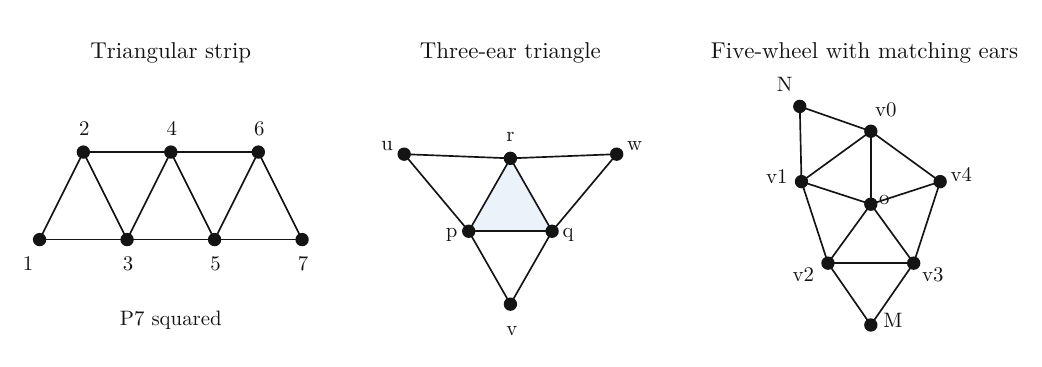}
\caption{Core-graph templates: the triangular strip $P_7^2$, the
three-ear triangle with central face $pqr$, and a five-wheel with ears
$N,M$ on disjoint rim edges. These drawings display the core graphs. The subsequent exclusions
concern positive midpoint realizations of their boundary completions.}
\label{fig:templates}
\end{figure}

Figure~\ref{fig:templates} shows three of the core templates.

\begin{lemma}\label{lem:doubleinterior}
The $I=2$ alternative is not terminal and cannot be the entire
two-connected $G$.
\end{lemma}
\begin{proof}
Use the labels in Lemma~\ref{lem:blocktemplates}. The eight star faces
are
\[
 opa,opb,oac_1,oc_1c_2,oc_2b,pbd_1,pd_1d_2,pd_2a.
\]
The saturated vertices $a,b$ have unique exterior neighbors $A,B$.
They cannot be cores: adjacency to two block vertices would place
them in the maximal block and violate its boundary degree bound.
The full remaining stars, with exterior thirds $n,m$ on $c_1c_2$
and $d_1d_2$, force
\[
 ac_1A,ad_2A,bc_2B,bd_1B,
 c_1c_2n,c_1An,c_2nB,d_1d_2m,d_1Bm,d_2mA.
\]
Here $A\ne B$, since otherwise flagness gives a third face on
$c_1c_2$. Also $n\ne m$: equality gives six distinct neighbors
$c_1,c_2,B,d_1,d_2,A$ and a closed link, impossible for either an
interior degree-at-most-five point or a boundary point. Other
identifications repeat an original link slot. Thus the eighteen faces
fill an original disk bounded by $(A,n,B,m)$, with all eight seed
vertices strictly inside. Each of $A,B$ has five known neighbors
and is noncore, so both are boundary. Each of $n,m$ has four known
neighbors and cannot be a four-center, whose closing link would
force the forbidden diagonal $AB$.

Suppose $m$ is boundary. If $n$ is boundary too, the outer cycle
is the physical rectangle. Otherwise a flat $m$ would leave $n$
as the sole outer point away from its physical side, violating the
two-sign rule at $n$. Normalize the true corner by
$m=0,A=(1,0),B=(0,1)$. If $n$ is a core, its full star forces
$n\in(0,1)^2$ exactly as in Lemma~\ref{lem:twobridges}; if boundary,
$n=(1,1)$. All seed vertices are in the open unit square.
The eight faces involving $m,d_1,d_2,p,a,b,A,B$ contradict
Lemma~\ref{lem:cornerobstruction}. Thus both $n,m$ are cores.
Lemma~\ref{lem:twobridges} supplies two distinct bridges. The eight
seed vertices and the two ear tips form a two-connected graph with
all seed stars complete and only those exits, so it is exactly the
maximal block and has two cuts.
\end{proof}

\subsection{Long triangular strips}
\begin{lemma}\label{lem:cornerfivefan}
Suppose a noncut five-star has cyclic link $p,q,U,k,h$, with
$Uk,kh$ physical and $k$ a true corner. Normalize
$k=0,U=(1,0),h=(0,1)$ and its center $v=(1/2,1/2)$.
Then its two other neighbors have the forced coordinates
\[
 p=(x,R),\quad q=(x+\tfrac12,R-\tfrac12),\quad
 R=\frac{3x}{4x-1},\qquad x>\tfrac14.
\]
\end{lemma}
\begin{proof}
Write $q=v+(x,y)$. Its equation with $U$ is
$(4x-1)y=1-x$; positivity of $q_x,q_y$ gives $x>1/4$.
For $u=p_x-1/2$ the other two equations reduce to
\[
 2(2x-1)u^2-(4x^2-6x+1)u-x(2x-1)=0.
\]
The discriminant is $(4x^2-2x+1)^2$. The candidates are
\[
 p_1=(x,3x/(4x-1)),\qquad
 p_2=(1/[2(1-2x)],3/[2(2x+1)]).
\]
For $x>1/2$ the second is not interior. For $1/4<x<1/2$,
\[
 \det(p_2-v,q-v)=
 \frac{2x(1-x)(4x^2-2x+1)}{(1-2x)(4x-1)(2x+1)}>0,
\]
whereas the actual link requires a negative determinant; $p_1$
has determinant $-(4x^2-2x+1)/(2(4x-1))<0$.
At $x=1/2$ the equation is linear and gives only $p_1$.
\end{proof}

\begin{lemma}\label{lem:boundarypoles}
Let a maximal core block have filled disk $P_n^2$, $n\ge5$, and
let $v_1$ be a noncut endpoint. If the sole exterior neighbor $U$
of $v_3$ and the exterior third $V$ on
$v_2v_4$ are boundary, the strip is impossible. The fifth neighbor
of $v_4$ may have any type.
\end{lemma}
\begin{proof}
The full stars give actual faces $v_1v_3U,v_3v_5U,v_2v_4V$.
The four distinct known neighbors of $v_4$ are $v_2,v_3,v_5,V$;
write $e$ for its remaining neighbor, without assuming its type.
The poles differ, since otherwise flagness gives a third face on
$v_3v_4$. The complete links at $v_2,v_1$ are
\[
 v_1,v_3,v_4,V,h,\qquad v_2,v_3,U,k,h.
\]
The vertices $h,k$ are noncores: a core $h$ would join the block
at $v_1,v_2$, and a core $k$ would make $v_1$ a cut or an extra
block neighbor. They differ from the displayed vertices; $h=U$
would prematurely close the link at $v_1$. Independence excludes
both being four-centers.

If $h$ is a four-center, its quad is $(v_1,v_2,V,k)$. Its physical
base admits the axial normalization
\[
 k=0,V=(1,0),\quad v_1=(t,1),v_2=(1-t,1),\quad
 h=(\tfrac12,\tfrac{1-4t}{2(1-t)}),\quad
 -\tfrac12<t<\tfrac14.
\]
The actual face $v_1v_2v_3$ gives $(v_3)_x=1/2$.
If $k$ is flat, $U=(2t,0)$ with $t<0$; four neighbors of $v_2$
are strictly left of $v_2$, leaving only $v_4$ for two right slots.
Otherwise
\[
 U=(0,2),\quad0<t<\tfrac14,\quad
 v_3=(\tfrac12,R),\quad R=\frac3{2(1-t)}\in(3/2,2).
\]
The two-sign rule at $v_2$ forces $(v_4)_x>1-t$. The two equations
with thirds $v_3,V$ reduce to $2X^2-3X+1=0$, so
$v_4=(1,2)$; the root $1/2$ is too small. The actual faces at
$v_3v_5$ then give
\[
 (v_5)_y=4-R,\qquad3(2-R)((v_5)_x-1/2)=0.
\]
Thus $v_5=(1/2,4-R)$. At $v_4$ its neighbors $v_2,v_3,v_5,V$
are weakly left, leaving only $e$ strictly right.

If $k$ is a four-center, its quad is $(v_1,U,l,h)$. The fourth
corner is boundary, by independence and connected-core separation
across the BB diagonal $Uh$. Normalize the corner quad as
$l=0,U=(0,1),h=(1,0),v_1=(1,1),k=(1/2,1/2)$.
Here $l\ne V$, since otherwise the physical edge $hV$ has two
different face-thirds. Thus $hV$ extends along the bottom and
\[
 v_2=(P,1/2),\quad V=(2P-1,0),\quad P>1,
 \qquad v_3=(1/2,3/(4P)).
\]
Every neighbor of $v_1$ now has ordinate at most one, a contradiction.

Finally $h,k$ are boundary. The two physical bases at $k$ cannot
be collinear. Use Lemma~\ref{lem:cornerfivefan}, putting
$D=4x-1$ and $H=2R-1=(2x+1)/D$.
The physical edge $hV$ is different from $hk$: $V=k$ would give
two face-thirds on that physical edge. If $h$ is flat, $V=(0,H)$,
$H>1$, hence $x<1$. The equations at $v_4$ give
\[
 Y=(X+1)/D,\quad
 2X^2+(1-6x)X+4x^2-2x=0,
\]
with candidates $(2x,H)$ and $(x-1/2,H/2)$. The latter is outside
if $x\le1/2$; if $x>1/2$, the full star at $v_3$ requires both
$v_4,v_5$ strictly to its right and excludes it. At $x=1/2$ the
same rule excludes even the first candidate. If $h$ is a corner,
the actual global height is one, $V=(2x,1)$, and $R<1$ forces
$x>1$. The equations instead give
\[
 Y=(8x+1-3X)/D,\quad2X^2-(6x+1)X+4x^2+2x=0.
\]
Besides $(2x,H)$ the root has ordinate $(10x-1)/(2D)>1$.
Thus in either case $v_4=(2x,H)$, $x\ne1/2$.
Now $v_3$ is the midpoint of $U,v_4$. For
$a=v_4-v_3=(x-1/2,H/2)$ and $z=v_5-v_3$, the two actual equations
give $B(a,z)=0$, $q(z)=-q(a)$, whence
$z=\pm(a_x,-a_y)$. The candidate $v_5=(2x,0)$ is boundary;
therefore $v_5=(1,H)$. In the flat case, four neighbors of $v_4$
have ordinate at most $H$, leaving only $e$ above. In the corner
case, $v_1,v_2,U,v_5$ are all strictly left of $v_3$, leaving only
$v_4$ right. Both contradict the two-sign rule.
\end{proof}

\begin{corollary}\label{cor:stripbound}
A terminal strip, including the whole two-connected case, has $n\le7$.
\end{corollary}
\begin{proof}
For $3\le i\le n-2$, the saturated vertex $v_i$ has a unique
noncore exterior neighbor $z_i$ and faces $v_iv_{i-2}z_i,v_iv_{i+2}z_i$.
It is noncore by maximality and Lemma~\ref{lem:cutslots}.
The unique exterior face on $v_iv_{i+2}$ gives $z_i=z_{i+2}$.
A pole shared by two such indices has four consecutive core neighbors
$v_{i-2},v_i,v_{i+2},v_{i+4}$; if it were a four-center, closing
its link would add the forbidden edge $v_{i-2}v_{i+4}$.
For $n\ge8$ both parity poles are consequently boundary. There is
a noncut endpoint by Lemma~\ref{lem:cutslots}, so
Lemma~\ref{lem:boundarypoles} applies.
\end{proof}

\begin{lemma}[Extreme endpoint of a four-core strip]\label{lem:extremeear}
In a terminal $P_4^2$, a noncut endpoint which is a coordinate extreme
has three boundary neighbors. Their two consecutive physical edges
meet at a true corner. With reflections and positive axis scalings,
its coordinates are those of Lemma~\ref{lem:cornerfivefan}, with
$x\ge1/2$.
\end{lemma}
\begin{proof}
Reflect to make the endpoint $v$ an abscissa minimum. Write its link
$p,q,a,b,c$, where $p,q$ are its two core neighbors, and $r$ is the
fourth block vertex, adjacent to both. The other three slots are
noncore. If $a$ is a four-center, its quad is $(v,q,d,b)$.
Its fourth corner cannot be $p$, since the original edge $vp$
would be a diagonal of that full four-star.
If $d$ were a core outside the block, $q$ would be a cut of block
degree three, contrary to Lemma~\ref{lem:cutslots}. If it were the
core $r$, the actual faces $qvp,qpr,qra,qav$ would already close
a four-cycle link at the original degree-five vertex $q$.
Thus such a fourth corner is boundary. The same applies at $c$.

If $b$ is a center, its quad is $(v,a,d,c)$. Independence makes
$a,c$ boundary, and connected-core separation across $ac$ makes
$d$ boundary. The positive four-star is a physical corner rectangle.
At a right corner the five neighbors of $v$ cannot include two
strictly left neighbors. At a left corner normalize
$d=0,a=(1,0),c=(0,1),v=(1,1)$. The actual face $vcp$ forces
$p_x=1/2$, contrary to extremality.

If $a,c$ are both centers, their fourth corners $d,e$ are boundary.
They are distinct, since equality fills an original triangle $pqd$
with the two four-stars and $vpq$, contradicting flagness.
The edges $db,be$ are physical. If $b$ is flat, axiality puts
$vp,vq$ on the same line, degenerating $vpq$. If $b$ is a corner,
they are perpendicular nonzero axial segments, so $E(v,p,q)\ne0$.

It remains to exclude a single end center, say $a$, with boundary
$b,c,d$. If $bd$ is vertical, then $vq$ is vertical and both are
block minima. Normalize this quad by
\[
 A=0,\quad B=(0,2),\quad
 z=((1-2t)/(2-t),1),\quad
 V=(1,t),\quad W=(1,2-t),\quad -1<t<1/2.
\]
Here $W=v$ is the degree-two companion and $V=q$ the degree-three
end; a reflection of the ordinate arranges this naming. The third
core on $VW$ is $p=(P,1)$, $P>1$, and $V$ has its other block
neighbor $r=(Q,Y)$. The two-sign rule at $V$ gives $Q>1$ and
$0<t<1/2$. The remaining neighbor $n=c$ of $W$ is boundary.
At $W$, its known abscissas are $0,z_x,1,P$, so $n_x>1$.
Consequently $Bn$ is physical horizontal and $n=(2,2)$.
The actual equations give
\[
 P=\frac3{2-t},\quad
 Y=\frac{t(Q-2)}{2Q-1},\qquad
 f(Q)=2Q^2+[1-P(4-3t)]Q+P(2-3t)-1=0.
\]
Since $f(1)=2(1-P)<0$ and $f(2)=0$, the only root above one is
$Q=2$, whence $Y=0$. The full link at $p$ contains
$r,V,W,n$ in this order and a fifth vertex $u$. The two actual
faces $pru,pnu$ imply
\[
 u_x=4-P,\qquad3(2-P)(u_y-1)=0.
\]
Here $P<2$, so $u_y=1$. At $r$ its four distinct neighbors
$A,V,p,u$ have ordinates $0,t,1,1$, and only its remaining slot
could lie below $r_y=0$, a contradiction.

If $bd$ is horizontal, positivity and $q_x\ge v_x$ give the axial
normalization
\[
 b=0,d=(1,0),\quad v=(t,1),q=(1-t,1),\quad
 a=(1/2,(1-4t)/(2(1-t))),\quad-1/2<t<1/4.
\]
If $b$ is flat, $c=(2t,0)$ and $t<0$, whereas $p_x=1/2$;
then only $c$ can be strictly left of $v$. Otherwise
$c=(0,2)$, $0<t<1/4$, and the actual face $vcp$ gives
$p=(1/2,R)$, $R=3/[2(1-t)]\in(3/2,2)$.
The full star at $q$ is $v,p,r,d,a$. Its two right slots force
$r_x>1-t$. Its two actual face equations reduce to
$2r_x^2-3r_x+1=0$, giving $r=(1,2)$.
Thus $c,d$ are flat points of the actual rectangle. At $p$ the
remaining neighbor $u$ is between $r,c$, giving faces $pru,pcu$.
It cannot be core, since edges $pu,ru$ would put it in the block.
If boundary, the physical edge $cu$ and $E(p,c,u)=0$ give
$u=(0,2R-2)$, strictly inside the already existing edge $bc$.
If a center, its quad $(p,r,k,c)$ has horizontal diagonal $rc$;
the four-star diagonal identity forces $pk$ horizontal as well,
contrary to strict convexity. This exhausts its types.

Therefore $a,b,c$ are all boundary. If $b$ were flat the two
axial bases would force $a=c$. A right corner leaves at most one
strictly left neighbor. At the left corner apply
Lemma~\ref{lem:cornerfivefan}; $p_x\ge v_x$ gives $x\ge1/2$.
The same proof, after interchanging coordinates and reflecting,
treats any of the four extreme labels.
\end{proof}

\subsection{The remaining short triangular strips}

We continue with a terminal core block, meaning a proper block with one
cut or the whole core graph when it is two-connected, after empty caps
have been removed. Write \(P_n^2\) for the strip with core vertices
\(1,\ldots,n\), edges \(ij\) exactly for \(1\le |i-j|\le2\), and
actual faces \((i,i+1,i+2)\). The numeral \(i\) also denotes its embedded
point. In face inventories below put
\[
 \mathcal T_n=\{(i,i+1,i+2):1\le i\le n-2\},\qquad
 \mathcal Q(z;a,b,c,d)=\{zab,zbc,zcd,zda\}.
\]
The second notation always means an actual complete four-star.

\begin{lemma}\label{lem:fourstrip}
A terminal block cannot be \(P_4^2\).
\end{lemma}
\begin{proof}
The cut, if any, is \(1\) or \(4\), by Lemma~\ref{lem:cutslots}.
First suppose the noncut endpoint \(1\) is a coordinate extreme.
Lemma~\ref{lem:extremeear} gives the all-boundary corner case.
Reflecting coordinates and using
Lemma~\ref{lem:cornerfivefan}, write
\[
 \begin{gathered}
 K=(0,0),\ A=(0,1),\ D=(1,0),\quad 1=(1/2,1/2),\\
 2=(x,R),\quad 3=(x+1/2,R-1/2),\quad
 R=\frac{3x}{4x-1},\quad x\ge\frac12.
 \end{gathered}
\]
The link at \(1\) is \((2,3,D,K,A)\). Since \(3\) is interior,
\(D\) is flat on the bottom side. Let \(B,C\) be the exterior thirds
of \(24,34\). The complete links are
\[
 2:(1,3,4,B,A),\qquad 3:(1,2,4,C,D),\qquad
 4:(2,3,C,z,B).
\]
Both \(B,C\) are noncore: a core attached to two block vertices would
belong to the same maximal block. They are distinct from one another
and from the opposite boundary labels, since otherwise flagness gives
a third face on \(23\), or closes a three-cycle link at \(4\).
No type is assumed for \(z\).

If both \(B,C\) are four-centers, their quads are \((2,4,z,A)\) and
\((3,D,z,4)\). A boundary \(z\) would be physically adjacent to both
\(A,D\). The bottom-flat point \(D\) forces \(z\) on the bottom,
so \(z=K\); but then \(AK\) has the two distinct incident faces
\(1AK,BAK\), impossible for a physical edge. If \(z\) is a core,
its five known neighbors \(4,B,A,D,C\) force its full link to close
\(A-B-4-C-D\), adding the edge \(AD\). This edge contains \(1\)
in its relative interior, impossible. Independence excludes a center
\(z\). At least one of \(B,C\) is therefore boundary.

Suppose \(B\) is boundary. If \(A\) is flat, its physical continuation
and face \(2AB\) give \(B=(0,2R-1)\) and \(x<1\). If \(A\)
is the upper-left corner, the global height is one, giving \(x>1\)
and \(B=(2x,1)\). Put \(d=4x-1\). Eliminating the two actual
face equations \(E(2,3,4)=E(2,B,4)=0\), for \(4=(X,Y)\), gives
\[
\begin{array}{c|c|c}
 &Y&\text{equation for }X\\ \hline
 A\text{ flat}&(X+1)/d&2X^2+(1-6x)X+4x^2-2x=0\\
 A\text{ corner}&(8x+1-3X)/d&2X^2-(6x+1)X+4x^2+2x=0.
\end{array}
\]
The first row gives \(4=(2x,2R-1)\) or
\((x-1/2,R-1/2)\). In the full star at \(3\), neighbors \(1,2,D\)
are weakly left, so both \(4,C\) must be strictly right. This excludes
the second root and also \(x=1/2\). The second row gives the same
first root or \((x+1/2,(10x-1)/(2d))\); the latter is above the
global height one. Thus \(x>1/2\) and \(4=(2x,2R-1)\).
If \(C\) is a center, its quad \((3,4,z,D)\) has \(3=(D+4)/2\),
contradicting strict convexity. If \(C\) is boundary, its physical
bottom edge with \(D\) gives \(C=(2x,0)\). The four known neighbors
\(2,3,B,C\) of \(4\) are then weakly left; its one remaining slot
cannot supply two strict right neighbors.

It remains that \(C\) is boundary and \(B\) a center. Again
\(C=(2x,0)\), \(x>1/2\). The equations
\(E(2,3,4)=E(3,C,4)=0\) factor to give
\(4=(2x,2R-1)\) or \((x,R-1)\). The full star at \(3\) excludes
the second, while the first has \(2=(A+4)/2\) in the full quad at
\(B\), again impossible. No noncut endpoint can be an extreme.

Such an endpoint must nevertheless exist. By the simultaneous-extreme
charge bound, a boundary core of disk degree three cannot attain two
coordinate-extreme labels. Thus \(2,3\) cannot account for all four
labels. In the whole-block case this already forces \(1\) or \(4\).
In the proper case, if only the cut \(4\) accounts for the remaining
labels, it must be a simultaneous extreme. Reflect to make it a double
minimum. The following elementary face-box observation applied to the
faces \(234,123\) then makes \(1\) a coordinate maximum.

For that observation translate a common double minimum \(u\) to zero
in two actual triangles \(rsu,rsv\), write \(r=(a,b),s=(c,d)\ge0\),
\(M=\max(a,c)>0,N=\max(b,d)>0\). The equation at \(u\) is
\(2ab+2cd=ad+bc\), and nondegeneracy implies
\((a+c)/M+(b+d)/N>2\); equality would put \(r,s\) on different
axes and violate this equation. For \(v=(X,Y)\), subtraction gives
\(2XY=(a+c)Y+(b+d)X\). Nondegeneracy and the minimum hypothesis give
\(X,Y>0\). If \(X\le M,Y\le N\), division by \(XY\) contradicts
the preceding strict inequality. Hence \(v\) is larger than both
\(r,s\) in at least one coordinate. This includes tied extrema and
completes the contradiction.
\end{proof}

\begin{lemma}\label{lem:fivetoseven}
A terminal block cannot be \(P_n^2\) for \(n=5,6,7\).
\end{lemma}
\begin{proof}
A saturated core \(i\) has its sole outside neighbor \(z_i\) and
actual faces \(i(i-2)z_i,i(i+2)z_i\). These poles are boundary or
four-centers, and \(z_i=z_{i+2}\) when both indices are saturated.
A shared four-center would require the forbidden edge
\((i-2)(i+4)\) to close its full link, so a shared pole is boundary.
We use the boundary-pole end lemma, Lemma~\ref{lem:boundarypoles},
including its version with an arbitrary fifth neighbor at vertex \(4\).

For \(n=7\), let \(U=z_3=z_5\), necessarily boundary, and \(V=z_4\).
Two boundary poles are excluded by that lemma. If \(V\) is a center,
its quad is \((2,4,6,A)\), where \(A\) is boundary. Indeed a core
fourth corner would belong to this maximal block, whereas its only
possible common neighbor of \(2,6\) is the already used corner \(4\).
Independence excludes a center corner. The full links force the faces
\[
 \mathcal T_7,\quad U13,U35,U57,\quad
 \mathcal Q(V;2,4,6,A),\quad A12,A67.
\]
Here \(U\ne A\): equality closes the link \(2-3-U-2\) at \(1\).
These fourteen actual faces have complete interior links and outer cycle
\((U,1,A,7)\); each core port has disk degree four. Lemma~\ref{lem:twobridges}
forces two distinct outgoing bridges, contrary to terminality.

For \(n=6\), put \(U=z_3,V=z_4\); they are distinct, since equality
gives a third flag face on \(34\). If both are centers, their quads
are \((1,3,5,A)\), \((2,4,6,B)\), with boundary fourth corners by
the same maximal-block argument. The full links at \(2,5\) force
\(12B,56A\). Here \(A\ne B\), since equality closes a four-link at
\(1\). The full links at \(1,6\) now force \(1AB,6AB\). A connected
core graph cannot lie on both sides of the proper boundary chord \(AB\),
so this is impossible.

The remaining six-core case, after reversal, has \(U\) boundary and
\(V\) a center with quad \((2,4,6,A)\). Let \(b\) be the exterior
third on \(56\). Again \(U\ne A\), since equality closes the same
three-cycle link at \(1\) as in the seven-core case. The links at
\(5,6\) are
\((3,4,6,b,U)\), \((4,5,b,A,V)\), forcing \(56b,5bU,6bA\).
A core \(b\) would join the maximal block. A center \(b\) would close
its known link path \(5-6-A-U\) with \(UA\); flagness would then close
a four-link at \(1\) using \(123,13U,12A,1UA\). Thus \(b\) is boundary.
The fourteen faces
\[
 \mathcal T_6,\quad U13,U35,\quad\mathcal Q(V;2,4,6,A),
 \quad A12,56b,5bU,6bA
\]
fill the simple cycle \((U,1,A,b)\). Its strict interior vertices are
\(2,3,4,5,6,V\), and port \(1\) has four disk neighbors. The physical
edges \(Ub,bA\) make \(b\) a true corner: if flat, the four known
neighbors of \(1\) lie strictly inward from it in the side-normal
coordinate, leaving one insufficient outward slot. Normalize
\(b=0,A=(1,0),U=(0,1)\). The same slot argument gives
\(1\in(0,1)^2\), hence all strict interior seed vertices lie there.
The eight actual faces
\[
 b6A,b5U,b56,456,\quad A6V,46V,U53,453
\]
contradict Lemma~\ref{lem:cornerobstruction}, with
\((m,d_1,d_2,p,a,b')=(b,5,6,4,V,3)\).

It remains to treat \(n=5\). Let \(U=z_3\), and let \(V,h,k\) be
the exterior thirds on \(24,12,45\). They are noncore, and the full
links at \(2,4\) are
\[
 (1,3,4,V,h),\qquad(2,3,5,k,V).                     \tag{a}
\]
If \(U\) is boundary, choose the noncut endpoint as \(1\). The
boundary-pole end lemma forces \(V\) to be a center, with quad
\((2,4,k,h)\); thus \(h,k\) are boundary. The link at \(1\) is
\((2,3,U,l,h)\), where \(l\) is noncore. If \(l\) is a center,
its private quad \((1,U,a,h)\) has a boundary corner \(a\), by the
opposite-core crosscut obstruction. Normalize its corner rectangle to
\(a=0,U=(0,1),h=(1,0),1=(1,1)\). The edge \(hk\) is physical
and distinct from \(ha\), for otherwise \(ha\) has two interior thirds.
Since \(h\) is flat, \(k=(K,0)\), \(K>1\). The axial quad at \(V\)
gives \(x_2+x_4=1+K\), and actual face \(234\) gives
\(x_3=(1+K)/2>1\). Face \(1U3\) instead gives \(x_3=1/2\).

Thus \(l\) is boundary and is a true corner. Normalize
\(l=0,U=(1,0),h=(0,1),1=(1/2,1/2)\). The corner-five-fan rule gives
\[
 2=(x,R),\quad3=(x+1/2,R-1/2),\quad
 R=\frac{3x}{4x-1},\quad x>1/4.
\]
Again \(k\ne l\), as the edge \(hl\) already has its incident face
at \(1\). If \(h\) is flat, \(k=(0,K)\), \(K>1\). The axial quad
and \(234\) force \(4=(x,R-1)\), so the directed opposite sides
\(24\) and \(kh\) point the same way, violating strict convexity.
If \(h\) is a corner, the global height is one, so \(R<1,x>1\).
Write \(k=(K,1)\). The same equations give \(4=(x+1,R)\),
\(K=2x+1\). Reflecting and scaling this axial quad gives bottom corners
\((0,0),(1,0)\) and upper corners \((t,1),(1-t,1)\), where
\(t=x/(2x+1)>1/3\). Its full four-star center would have ordinate
\((1-4t)/(2(1-t))<0\), by substituting its midpoint abscissa in an
actual face equation. It cannot lie inside the quad. This excludes
boundary \(U\).

Now \(U\) is a center with quad \((1,3,5,A)\), where \(A\) is
boundary by maximality. The end links are
\((2,3,U,A,h)\), \((3,4,k,A,U)\). In particular \(A\ne h,k\).
If \(V\) is a center, independence makes \(h,k\) boundary. The faces
\[
 \mathcal T_5,\quad\mathcal Q(U;1,3,5,A),\quad
 \mathcal Q(V;2,4,k,h),\quad12h,1Ah,45k,5kA
\]
fill the original triangle \(Ahk\), contradicting flagness.
Thus \(V\) is boundary. If \(h\) is a center, its quad is
\((1,2,V,A)\). If \(k\) is also a center, the two centers provide
actual faces on opposite sides of the same boundary chord \(AV\),
separating the connected cores. If \(k\) is boundary, the faces
\[
 \mathcal T_5,\quad\mathcal Q(U;1,3,5,A),\quad
 \mathcal Q(h;1,2,V,A),\quad24V,45k,5kA,4kV
\]
fill the forbidden original triangle \(AVk\). Reversal covers the
case with \(k\) a center. Consequently \(V,h,k\) are all boundary.
Here \(A\ne V\), since \(12A\) would be a third face on \(12\),
and \(h\ne k\), since equality closes a six-cycle link at a boundary
vertex. The fourteen faces
\[
 \mathcal T_5,\quad\mathcal Q(U;1,3,5,A),\quad24V,
 \quad12h,1Ah,2Vh,45k,5kA,4kV
\]
fill the physical four-cycle \((A,h,V,k)\), which is the entire
rectangular boundary. All seed vertices are strictly inside. At corner
\(V\), the eight faces
\[
 V2h,V4k,V24,324,\quad h21,321,k45,345
\]
contradict Lemma~\ref{lem:cornerobstruction}, with
\((m,d_1,d_2,p,a,b')=(V,4,2,3,1,5)\).
All inventories above are actual faces. Their complete interior links
and boundary cancellation identify the indicated filled disks; no
artificial half of a quadrilateral was used.
\end{proof}

\subsection{The three-ear triangle}

The exclusion is organized by the number of center poles. The first
lemma leaves only zero or one. An ear-translation lemma is followed by
a separate exclusion of each remaining case. No all-center geometric
case is left for a coordinate argument.

Assume that a terminal block has the three-ear form.
Let the central triangle be \(pqr\), with ears \(v,w,u\) on
\(pq,qr,rp\). The degree-four central cores have unique outside
noncore neighbors \(A,B,C\), respectively. Their full links are
\[
 p:(q,r,u,A,v),\quad q:(r,p,v,B,w),\quad
 r:(p,q,w,C,u).
\]
The ear links are
\[
 v:(p,q,B,h_v,A),\quad w:(q,r,C,h_w,B),\quad
 u:(r,p,A,h_u,C).                                   \tag{b}
\]

\begin{lemma}[Only boundary poles or one center]\label{lem:sunpoletypes}
In the displayed three-ear configuration, the poles are distinct and
at most one is a four-center. If $A$ is the only center, its common
port $h=h_v=h_u$ is boundary.
\end{lemma}
\begin{proof}
The poles \(A,B,C\) are pairwise distinct: otherwise a central edge
has a third flag face. No port \(h_i\) equals a pole. Equality with
an adjacent pole repeats a link vertex; equality with the opposite pole
would add an edge between the two poles in another ear's link, forcing
an absent flag face. The known link path at \(A\) is
\(h_v-v-p-u-h_u\). If its ends are distinct, \(A\) is boundary;
if they coincide, they close a four-cycle and \(A\) is a center.
Thus
\[
 A\text{ center}\iff h_v=h_u,\quad
 B\text{ center}\iff h_v=h_w,\quad
 C\text{ center}\iff h_w=h_u.
\]
Two center poles would make all three ports equal and close a six-cycle
link there, impossible. If \(A\) is the single center, its common port
\(h=h_v=h_u\) is boundary: independence excludes a center, and a core
attached to both \(u,v\) would belong to the maximal block.
\end{proof}

\begin{lemma}[Boundary-pole ear translation]\label{lem:suneartranslation}
For a noncut ear with two boundary poles, the remaining neighbor is a
true boundary corner. With the adjacent core labels as below, the
translation identity \eqref{eq:suntranslation} holds.
\end{lemma}
\begin{proof}
If a noncut
ear \(e\) has adjacent core neighbors \(a,b\) and boundary poles
\(L,R\), its remaining neighbor \(h\) is noncore. A center \(h\)
would have quad \((e,L,d,R)\). Its fourth corner \(d\) is boundary,
since opposite cores across the proper boundary chord \(LR\) violate
connectedness. The physical corner quad then gives
\(d=0,L=(\alpha,0),R=(0,\beta),e=(\alpha,\beta)\), with
\(h=(\alpha/2,\beta/2)\). Face \(ebR\) gives \(b_x=\alpha/2\).
Neighbors \(b,R,h\) are strictly left of \(e\), and \(L\) is equal,
leaving only one possible right slot. Thus \(h\) is boundary and a
true corner. Lemma~\ref{lem:cornerfivefan} then gives, with the core
labels chosen adjacent to their corresponding poles,
\begin{equation}\label{eq:suntranslation}
 b-a=(R-L)/2.
\end{equation}
\end{proof}

\begin{lemma}[Three boundary poles]\label{lem:sunboundary}
A terminal three-ear triangle cannot have all its poles on the boundary.
\end{lemma}
\begin{proof}
Two ears are noncut;
rotate labels to make them \(v,u\), sharing pole \(A\). Their ports
are distinct true corners, and \eqref{eq:suntranslation} gives
\[
 q-p=(B-A)/2,\qquad r-p=(C-A)/2,
 \qquad E(A,B,C)=0.
\]
\paragraph{The common pole is a corner.}
If \(A\) is a true corner, normalize the whole rectangle to
\(A=0,h_v=(1,0),h_u=(0,1)\), with \(B=(1,b),C=(c,1)\),
\(0<b,c\le1\). Then \(v=(1/2,b/2),u=(c/2,1/2)\).
For \(p=(X,Y)\), faces \(Apv,Apu\) give
\[
 (4X-1)Y=b(X-1),\qquad(4Y-1)X=c(Y-1).
\]
Since \(p\) is interior, these imply \(X,Y<1/4\), then \(Y>b,X>c\).
Consequently \(b,c<1/4\), contradicting
\(E(A,B,C)=2b+2c-1-bc=0\).

\paragraph{The common pole is flat.}
If \(A\) is flat, use the global height one and a local horizontal
scale to write
\[
 h_v=0,\ h_u=(0,1),\ A=(0,a),\ B=(1,0),\ C=(k,1),
 \qquad0<a<1.
\]
Put \(\beta=B-A,\gamma=C-A,s=p-A,t=w-p\).
The translation identities and the central face give
\(k=(1+a)/(2-a)\), \(q=p+\beta/2,r=p+\gamma/2\).
The actual faces \(Apv,Apu\) say
\[
 B(s,\beta)=4q(s)+q(\beta),\qquad
 B(s,\gamma)=4q(s)+q(\gamma).
\]
Expanding the actual outer faces \(wqB,wrC\) gives respectively
\[
 -4q(s)+2q(t)-B(t,3\beta/2-s)=0,\qquad
 -4q(s)+2q(t)-B(t,3\gamma/2-s)=0.
\]
Subtract these, and subtract the faces \(wqr,pqr\), to obtain
\(B(t,\beta-\gamma)=0\) and \(4q(t)=B(t,\beta+\gamma)\).
Writing \(\delta=1-k,T=t_y\), the first gives \(t_x=\delta T\).
If \(T=0\), then \(w=p\); if \(\delta=0\), the second equation
forces \(T=0\). Otherwise
\[
 T=\frac{1-a\delta}{2\delta}
  =\frac{1-a+a^2}{1-2a}.
\]
For \(a<1/2\) this is greater than one; for \(a>1/2\) it is less
than minus one. Both contradict \(0<p_y,w_y<1\). The zero-denominator
case was already excluded.
\end{proof}

\begin{lemma}[One center pole]\label{lem:sunonecenter}
A terminal three-ear triangle cannot have exactly one center pole.
\end{lemma}
\begin{proof}
Let \(A\) be the single center, with common boundary port
\(h=h_v=h_u\). All the faces in the displayed full stars, except
\(wBh_w,wCh_w\), form a fourteen-face disk with boundary
\((h,B,w,C)\), strict interior \(p,q,r,v,u,A\), and disk degree four
at \(w\). Its boundary edges \(hB,hC\) are physical. If \(h\) were
flat, the four known neighbors of \(w\) would be strictly inward,
leaving just one possible outward neighbor. Thus \(h\) is a true
corner. Normalize \(h=0,B=(1,0),C=(0,1)\). If either coordinate of
\(w\) were at least one, the same four known neighbors would be weakly
smaller, leaving one insufficient larger slot. Hence \(w\), and then
all strict interior seed vertices, lie in \((0,1)^2\).

\paragraph{Coordinate identities and strict bounds.}
Write \(v=(1/2,y),u=(x,1/2),A=(s,t),p=(P,T)\).
The actual faces \(hvA,huA\) give
\[
 0<s,t<1/4,\qquad
 x=\frac{s(1/2-2t)}{1-t},\qquad
 y=\frac{t(1/2-2s)}{1-s}.
\]
The full four-star identity \(B(p,u-v)=0\) gives
\[
 \frac{T}{P}=\frac{1/2-y}{1/2-x}=\frac{1-t}{1-s},
\]
where the second equality uses the preceding expressions for \(x,y\).
Subtracting
\(E(h,v,A)\) from \(E(p,v,A)\) then yields
\[
 N_0=\tfrac12+s+t-4st,\quad
 P=\frac{N_0}{2(1-t)},\quad T=\frac{N_0}{2(1-s)},
 \qquad \tfrac14<P,T<\tfrac12.
\]
For example \(N_0-(1-t)/2=s+3t/2-4st>0\), and
\((1-t)-N_0=(1/2-s)(1-4t)>0\); the other bounds are symmetric.
\paragraph{The remaining central face.}
For \(q=(X,Y)\), the actual face \(Bvq\) has equation
\[
 (4X-3)Y=2Xy,
\]
so \(X>3/4\); symmetrically \(r_y>3/4\). Set
\(d=1/2-P,f=X-1/2\), so \(0<d<1/4<f\). Face \(pvq\), translated
by \(v\), gives
\[
 (2f+d)(Y-y)=(f+2d)(T-y),\qquad
 T-y=\frac{1/2+s}{2(1-s)}>0.
\]
The ratio on the right is strictly less than one, hence \(Y<T\).
The symmetric argument gives \(r_x<P\). Thus \(q-p\) has signs
\((+,-)\) and \(r-p\) signs \((-,+)\). Their quadratic values are
negative and their bilinear pairing positive. The actual central face
would have
\[
 E(p,q,r)=2q(q-p)+2q(r-p)-B(q-p,r-p)<0,
\]
a contradiction. This excludes the single-center case.
\end{proof}

\begin{lemma}[Three-ear triangle exclusion]\label{lem:sun}
A terminal block cannot consist of a triangle with one triangular ear
on each of its edges.
\end{lemma}
\begin{proof}
Lemma~\ref{lem:sunpoletypes} leaves only three boundary poles or one
center pole. Lemmas~\ref{lem:sunboundary} and~\ref{lem:sunonecenter}
exclude these alternatives.
\end{proof}

\subsection{Wheel annuli and their last central faces}
\begin{lemma}\label{lem:wheelnoF}
In the one-interior-core template every exterior rim third is
boundary or a core, never a four-center.
\end{lemma}
\begin{proof}
Let $o$ have cyclic rim $v_0,\ldots,v_4$, and let $y_i$ be the
exterior third on $v_iv_{i+1}$. The full rim stars give all fifteen
actual faces
\[
 ov_iv_{i+1},\qquad v_iv_{i+1}y_i,\qquad v_iy_{i-1}y_i.
\]
The $y_i$ are distinct: adjacent equality repeats a slot and a
distance-two equality gives a third face on the intervening rim edge.
None is an old rim point, since a rim chord is incompatible with
the full link at $o$. Thus these faces form a disk with outer cycle
$(y_0,\ldots,y_4)$. A core $y_i$ belongs to the maximal block and
is one of its matching ear tips.

A four-center $y_i$ closes its four known neighbors by the actual
exterior face $y_{i-1}y_iy_{i+1}$. Two such centers are nonadjacent
by independence; adding their two exterior triangles leaves a
separating original three-cycle. If a four-center $y_3$ is adjacent
to a core $y_2$, its face $y_2y_3y_4$ gives $y_2$ five neighbors
$v_2,v_3,y_1,y_3,y_4$; its full link forces $y_2y_4y_1$.
Again the remaining boundary $y_0y_1y_4$ separates the seed wheel.
Therefore a possible single four-center has two boundary neighbors.
Its BB closing edge is physical, and it forms a Q wall whose two
opposite rim cores both have core degree three.

Here is the full local contradiction for that wall. Normalize its
physical base $A=(0,0),B=(0,2)$ and opposite cores
$v=(1,t),w=(1,2-t)$, initially $-1<t<1/2$ by the four-star rule.
Its center $z_0$ has ordinate one and abscissa between zero and one.
The other actual faces are $vwp,vpq,vqA,wpr,wrB$, with $p=o$.
At $v,w$ the two-sign rule gives $p_x,q_x,r_x>1$; the known
four ordinates at $v$ exclude $t\le0$. Thus $p=(P,1)$, $P>1$,
and $0<t<1/2$. The equations at $q=(X,Y)$ are
\[
 Y=\frac{t(X-2)}{2X-1},\qquad
 f(X)=2X^2+[1-P(4-3t)]X+P(2-3t)-1=0.
\]
Since $f(1)=2(1-P)<0$, there is one root $Q>1$; the reflected
equations force $q=(Q,Y),r=(Q,2-Y)$. The fifth neighbor $z$ of
$p$ is the remaining rim core. Subtracting $pqz,prz$ gives
$z=(2Q-P,1)$; $Q=P$ would give a degenerate triangle, and $Q<P$
puts every neighbor of $p$ to its left. Hence $Q>P$.
At $q$, the three neighbors $v,p,z$ have ordinate greater than $Y$;
if $Y\le0$, its fourth neighbor $A$ is also weakly above, impossible.
Consequently $Q>2$ and $0<Y<t$.

Let $n,m$ be the fifth neighbors at $q,r$. The full link at the
core $z=(Z,1)$, $Z=2Q-P$, is $q,p,r,m,n$; in particular $nm$
is an original edge and $n\ne m$. The two-sign rule gives
$n_x,m_x>Q$. For $n=(X,T)$ its actual faces with $A,z$ give
\[
 T=\frac{Y(X-2Q)}{2X-Q},\qquad X+ZT=Z(2-Y)-Q,
\]
and therefore
\[
 F(X)=2X^2+[Q-Z(4-3Y)]X+QZ(2-3Y)-Q^2=0.
\]
Since $F(Q)=2Q(Q-Z)<0$, the allowed root $N>Q$ is unique;
reflection gives $n=(N,T),m=(N,2-T)$. At $z$, its three neighbors
$p,q,r$ are left, so $N>Z$. Now $n$ has known neighbors $A,q,z,m$
with abscissas $0,Q,Z,N$. If it were interior, at most one remaining
slot could be strictly right. Thus $n,m$ are boundary. Their physical
edges $An,Bm$ force $T=0$, so $N=2Q$. It follows that
\[
 0=F(2Q)=3Q[3Q-Z(2-Y)],\qquad Z>3Q/2.
\]
But $f(Q)=0$ gives
\[
 P=\frac{2Q^2+Q-1}{(4-3t)Q-2+3t}>Q/2,
\]
since the denominator is positive and twice the numerator minus
$Q$ times that denominator is $4Q-2+3tQ(Q-1)>0$.
Thus $Z=2Q-P<3Q/2$, the required contradiction.
\end{proof}

Use the completed outer cycle and fifteen actual faces in the
proof of Lemma~\ref{lem:wheelnoF}. Its noncore vertices are now all boundary.

We separate the zero-ear, one-ear and two-ear cases. In the last case,
one normalization supplies two distinct lemmas according as the physical
boundary continues straight or turns at a corner. The two final sign
arguments use different center bounds and are kept separate.

\begin{lemma}[No ears]\label{lem:wheelzero}
The wheel template without ears is neither proper terminal nor the
entire two-connected core graph.
\end{lemma}
\begin{proof}
With zero ears, the five physical outer edges form a rectangle with
one flat vertex. A corner not adjacent to that flat point has its
two adjacent outer edges equal to a full width and height. All seed
vertices lie in that open box, and its eight wheel faces contradict
Lemma~\ref{lem:cornerobstruction}.
\end{proof}

\begin{lemma}[One ear]\label{lem:wheelone}
The wheel template with one ear is neither proper terminal nor the
entire two-connected core graph.
\end{lemma}
\begin{proof}
Write the outer cycle $(N,A,B,C,D)$ with only $N$
a core. The physical chain $A-B-C-D$ cannot have both middle points
flat: $N$ would be its sole point off one side, with only one possible
outward neighbor. Nor can both be corners. In that case write
$B=0,C=(1,0),A=(0,a),D=(1,d)$. The full star of $N$ gives
$N_y<\max(a,d)$; at the higher of the two ends the adjacent true
corner and all seed vertices satisfy the unit-box hypotheses of
Lemma~\ref{lem:cornerobstruction}. Thus, up to reflection,
\[
 C=0,\ B=(0,1),\ A=(0,H),\ D=(1,0),\ H>1,
 \quad N=(n,m)\in(0,1)\times(0,H).
\]
The last box follows from the two-sign rule at $N$, allowing tied
outer maxima. Put the center $p=(P,T)$ and cyclic rim
\[
 z=(Z,W),\ u=(U,V),\ (a,b),\ s=(x,1/2),\ t=(1/2,y).
\]
The boundary faces give
\[
 x+y-xy=1/4,\quad0<x,y<1/4,\quad
 a=\frac{x(H+1)}{2H-1},\quad b=\frac{H+1}{2}.
\]
The actual face $Dtz$ gives $(4Z-3)W=2Zy$, so $Z>3/4$ and
$W>2y$. Apply Lemma~\ref{lem:quadrantdichotomy} at $D$, reflecting
the first coordinate: across $Dz$ the opposite thirds are $t,N$.
Since $1/2>2(1-Z)$, the former has negative determinant and the
latter positive; hence $1-n<(1-Z)/2$ and $m>2W$.
At $A$ reflect the ordinate. The coordinates of $(a,b),B$ are
$(a,e),(0,2e)$, $e=(H-1)/2>0$. Across the edge to $(a,b)$ the
opposite third $u$ must have negative determinant, giving
$U>2a$, $H-V<e/2$. Across $Au$ the opposite thirds are $(a,b),N$;
the former has positive determinant, so $n>2U$.
Consequently
\[
 U<1/2,\quad a<1/4,\quad W<H/2,\quad
 V>(3H+1)/4.
\]
The center's full five-star gives $P<1/2$, since otherwise only
$z$ could be strictly right. Its actual faces $Cst,pst$ give
\[
 T=\frac{P(1/2+y)}{2P-x-1/2},\quad
 \tfrac14<P<\tfrac12,\quad T>y.
\]
The actual face $ptz$ gives
\[
 W-T=(y-T)\frac{Z+P-1}{2Z-P-1/2}<0.
\]
All strict factors have the displayed signs. Since $u$ is the
unique highest rim point, the center inside its rim has $V>T>W$.
For $\beta=V-T>0$, $\delta=T-W>0$, the remaining actual face has
\[
 E(p,u,z)=\beta(2U-P-Z)+\delta(U+P-2Z)<0,
\]
because $U<1/2$, $P>1/4$, $Z>3/4$, and $P<1/2$.
This excludes one ear without any assumption about its fifth neighbor.
\end{proof}

\begin{lemma}[A boundary cap and two continuations]\label{lem:wheeltwocap}
A terminal two-ear wheel, including the whole two-connected case,
has the boundary-corner normalization below. Its only remaining
physical continuations are $B=(L,0)$, $L>1$, and $B=(1,L)$, $L>0$.
\end{lemma}
\begin{proof}
Label the outer cycle $(N,A,M,B,C)$, with $A,B,C$
boundary. Choose a noncut ear $N$; such a choice exists both for a
proper terminal block and for the entire two-connected case. Its
fifth neighbor $h$ is noncore. If it is a four-center, its quad is
$(N,C,l,A)$; the fourth corner is boundary by independence and
connected-core separation. It differs from $B$, since otherwise
the physical edge $BC$ has two actual thirds. The corner-quad rule
gives $l=0,C=(1,0),A=(0,1),N=(1,1),h=(1/2,1/2)$.
Then $C$ is flat and $B=(L,0)$, $L>1$. Both rim neighbors of $N$
must have ordinate $>1$, forcing $M_y>1$. At $M$ its two inner
neighbors and $A,B$ are now all lower, leaving only one upper slot.
Thus $h$ is boundary. Its physical edges to $C,A$ meet at a true
corner and give
\[
 h=0,\quad A=(0,1),\quad C=(1,0),\quad N=(1/2,1/2).
\]
Write the cyclic rim $q,p,u,r,s$ and center $o=(K,H)$.
Lemma~\ref{lem:cornerfivefan}, with cyclic link $p,q,C,h,A$, yields
\[
 p=(x,R),\quad q=(x+1/2,R-1/2),\quad
 R=3x/(4x-1),\quad x>1/4.
\]
The remaining physical edge $CB$ either continues flat, $B=(L,0)$,
$L>1$, or turns, $B=(1,L)$, $L>0$. Write
$u=(U,V),r=(Q,Y),s=(P,T)$ and $M=(n,m)$.
The full star at $M$ bounds both its coordinates by the coordinate
maxima of $A,B,C$: its two seed neighbors are strictly inside the
outer cycle, so a maximum at $M$ would leave one outward slot.
The same box contains the seed wheel.
\end{proof}

For the next two lemmas use the labels and actual faces of
Lemma~\ref{lem:wheeltwocap}: the rim is $(q,p,u,r,s)$, the center
is $o=(K,H)$, and $u=(U,V),r=(Q,Y),s=(P,T)$.

\begin{lemma}[Two ears: flat continuation]\label{lem:wheeltwoflat}
The flat continuation in Lemma~\ref{lem:wheeltwocap} is impossible.
\end{lemma}
\begin{proof}
In the flat case this box is $(0,L)\times(0,1)$, so $x>1$ and
$3/4<R<1$, while $P=(L+1)/2$. The following consequences of
Lemma~\ref{lem:quadrantdichotomy} record every actual edge and its
opposite thirds:
\[
\begin{array}{c|c|c}
\text{anchor and reflected coordinate}&\text{edge; opposite thirds}&
\text{consequence}\\ \hline
A,\ 1-y&Ap;\ N,u&U>2x,\quad V>(1+R)/2\\
A,\ 1-y&Au;\ p,M&n>2U,\quad U<L/2\\
B,\ L-x&Bs;\ C,r&Q>(3L+1)/4,\quad Y>2T\\
B,\ L-x&Br;\ s,M&m>2Y,\quad 0<T<1/4,\ Y<1/2
\end{array}
\]
For the first test $N_x=1/2<x/2$ gives the positive determinant.
For the third it is $-(L-1)T<0$.
The remaining signs follow from the preceding strict inequalities.
\paragraph{Locating the center.}
The full star of $q$ forces $K>x+1/2$. Its actual central face
with $p$ gives, for $a=K-x-1/2>0$,
$(4a+1)(H-R+1/2)=a+1$, hence $H>1/2$.
The center then needs both $p,u$ above it, so $H<R<V$.
Its right-slot rule gives $K<P$, since $p,q,u$ are left of $P$.
Subtract the actual equations $Cqs,oqs$ and put $X=K-1>0$. Then
\[
 2XH=(x-1/2+P-1)H+(R-1/2+T)X,
 \qquad K>L/4+x/2+1/2.
\]
\paragraph{The final sign contradiction.}
Thus for $\beta=V-H>0$, $\delta=H-Y>0$,
\[
 E(o,u,r)=\beta(2U-K-Q)+\delta(U+K-2Q)<0;
\]
the two coefficients are respectively less than $-x/2-3/4$ and
$-L/2$. This contradicts its actual face equation.
\end{proof}

\begin{lemma}[Two ears: corner continuation]\label{lem:wheeltwocorner}
The corner continuation in Lemma~\ref{lem:wheeltwocap} is impossible.
\end{lemma}
\begin{proof}
In the corner case the true global width is one, so
$1/4<x<1/2$, $R>3/2$. The outer maximum bound forces $L>R$;
all seed points and $M$ lie in $(0,1)\times(0,L)$, and $T=L/2$.
The same face tests, now in the following charts, give
\[
\begin{array}{c|c|c}
A,\ y-1&Ap;\ N,u&0<U<x/2,\quad V>2R-1\\
A,\ y-1&Au;\ p,M&n<U/2,\quad m>2V-1,\ V<(L+1)/2\\
C,\ 1-x&Cq;\ N,s&P>(3+2x)/4,\quad L>4R-2>4\\
B,\ (1-x,L-y)&Bs;\ C,r&Q<2P-1,\quad Y>3L/4
\end{array}
\]
At $A$, $\det((x,R-1),(1/2,-1/2))<0$; the positive-determinant
extension of Lemma~\ref{lem:quadrantdichotomy} proves $V>1$
rather than assuming it. At $C$ the determinant with $N$ is
$(1-x-R)/2<0$; at $B$ the determinant with $C$ is positive.
\paragraph{Locating the center and the high rim.}
The full star of $p$ forces $H>R$; $opq$ gives
\[
 K-x=\frac{H-R+1}{4(H-R)+1}>1/4.
\]
The full center star further gives $K<P$ and $H<(L+1)/2$:
otherwise only $r$ can occupy the missing right or upper slot.
Put $k=1-K$, $a=1/2-x$, $d=1-P$, $b=R-1/2$, all positive.
Subtracting $Cqs$ from $oqs$ gives
\[
 2kH=(a+d)H+(b+L/2)k,\qquad H>L/4+R/2-1/4>L/4.
\]
The actual face $ors$ now yields
\[
 Q-K=\frac{(P-K)(Y+H-L)}{2Y-H-L/2}>0;
\]
its denominator is $>(L-1)/2$ and its numerator is positive.
Independently subtracting $Apu$ from $opu$ gives
\[
 2K(H-1)=(x+U)(H-1)+(R+V-2)K,
 \qquad H>(R+V)/2.
\]
\paragraph{The final sign contradiction.}
With $\alpha=K-U>0$, $\gamma=Q-K>0$, the unused central face is
\[
 E(o,u,r)=\alpha(Y+H-2V)+\gamma(2Y-V-H)>0.
\]
Indeed
$2V-H-Y<3/4-R/2<0$ and $V+H-2Y<1-L/2<0$ by the bounds just
proved. No sign for $V-H$ was assumed. This final contradiction
excludes the corner case.
\end{proof}

\begin{lemma}[Wheel exclusion]\label{lem:terminalwheels}
No wheel-with-matching-ears template is a proper terminal block or
the entire two-connected $G$.
\end{lemma}
\begin{proof}
Lemma~\ref{lem:blocktemplates} leaves zero, one or two matching ears.
Lemmas~\ref{lem:wheelzero} and~\ref{lem:wheelone} handle the first
two cases. Lemma~\ref{lem:wheeltwocap} reduces the last to the two
continuations excluded by Lemmas~\ref{lem:wheeltwoflat}
and~\ref{lem:wheeltwocorner}.
\end{proof}

\begin{proposition}[Terminal-block availability]\label{prop:availability}
If $G$ is two-connected then $G=K_3$. If $|G|>3$, then $G$ has
either a degree-one core or a proper terminal $K_3$.
\end{proposition}
\begin{proof}
For any nontrivial terminal block, including the whole two-connected
case, Lemmas~\ref{lem:noquadterminal}, \ref{lem:doubleinterior},
\ref{lem:blocktemplates}, \ref{cor:stripbound}, \ref{lem:fourstrip},
\ref{lem:fivetoseven}, \ref{lem:sun}
and~\ref{lem:terminalwheels} leave only $K_3$.
If $|G|>3$, the whole graph is therefore not two-connected. A leaf
of its block-cut tree is either a bridge, whose noncut endpoint has
core degree one, or a proper terminal nontrivial block, now $K_3$.
Lemma~\ref{lem:cutslots} gives the latter exactly two noncut triangle
vertices and one cut with its unique outgoing bridge. These are the
full geometric configurations needed for the subsequent R/D moves.
\end{proof}

\section{From geometric reduction to a combinatorial presentation}

\begin{theorem}[Finite geometric reduction]\label{thm:reduction}
Suppose a positive midpoint rectangle has interior degrees at most five
and a nonempty connected five-core graph. After legitimate empty-cap removals,
it admits a finite sequence of the R and D reductions, ending with
zero, two or three cores. Every intermediate realization is positive,
has interior degrees at most five, and has connected or empty core graph.
Each removed R or D rectangle meets its remainder along a complete
single-edge side. The reduction is the inverse of a flip between two
four-valent apices, and retains the original vertex labels.
\end{theorem}
\begin{proof}
When there are more than three cores, Proposition~\ref{prop:availability}
gives a core leaf or a proper terminal triangle. Apply
Proposition~\ref{prop:leafreduction} or
Proposition~\ref{prop:pendantreduction}, respectively. These remove two
or four cores and have precisely the stated preservation properties.
Empty caps are the full rectangular strips of Lemma~\ref{lem:cap}.
The core count strictly decreases, so the process terminates; the
remaining cases are Proposition~\ref{prop:terminalbases}.
\end{proof}

\begin{lemma}[A fixed direction in reverse assembly]\label{lem:axislock}
In the reverse sequence, all nonempty rectangular attachments have the
same direction, including the restored empty caps.
\end{lemma}
\begin{proof}
Each attachment uses a whole side consisting of one edge. Its two
endpoints become strict flat boundary vertices on the two perpendicular
sides of the larger rectangle. Neither perpendicular side is henceforth
a single edge. An interface flip removes the joining interior edge,
not either boundary vertex. Further parallel attachments retain these
flat vertices. Therefore no subsequent attachment can be perpendicular.
This argument reverses the complete sequence, including all empty-cap
removals; it does not assume that a forward removal preserves earlier
flat vertices. Subdivided sides of the terminal rectangle are already
ineligible for an attachment.
\end{proof}

\begin{figure}[htbp]
\centering
\includegraphics[width=\linewidth]{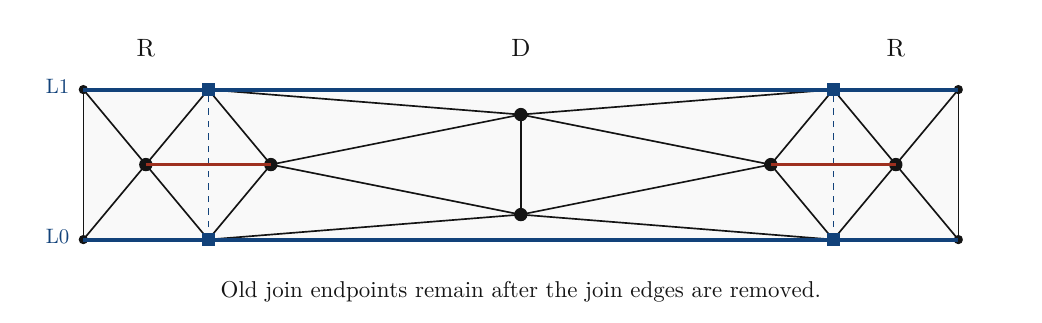}
\caption{An R--D--R presentation with both admissible joins flipped.
The blue outer paths are $L_0,L_1$; squares mark the original old join
endpoints. Dashed segments are the deleted join edges, while the red
segments are their replacements. All six interior vertices now have
degree five. The original old/new labels persist after each flip.}
\label{fig:serial}
\end{figure}

Figure~\ref{fig:serial} records the original endpoints of erased joins.

\begin{proposition}[Presentations supporting the vertex conditions]
\label{prop:presentations}
For a relative disk with interior degrees at most five and empty or
connected five-core graph, positive common weights satisfying \(V\)
exist exactly for the following marked presentations:
\begin{enumerate}
\item a serial word of R and D rectangles, with designated joins flipped
between their facing four-centers, each center used at most once;
\item a T rectangle with a possibly empty unflipped R strip along
\(DA\) or \(CD\), but not both.
\end{enumerate}
A D piece may use either opposite side pair as its ends, but a flipped
end must have a four-center apex. In either presentation \(L\) is a pair
of opposite complete outer side paths. The resulting core graph must
satisfy the stated connectivity hypothesis.
\end{proposition}
\begin{proof}
Starting from weights, Proposition~\ref{prop:rectangle} gives a positive
rectangle. For a nonempty connected core apply Theorem~\ref{thm:reduction}
and reverse it; an empty core is directly covered by
Proposition~\ref{prop:terminalbases}. The terminal
core-free rectangle is an R strip. A terminal two-core piece is D or F,
where F is two R pieces with their common join flipped. A terminal
three-core piece is T, with its allowed empty strip. Every restored
nonempty cap is R or D, and Lemma~\ref{lem:axislock} makes the complete
assembly serial. The facing apices have degree four before a flip and
five afterwards; using either twice would create degree six.
An unflipped D end with a degree-five apex is permitted, since ordinary
gluing does not change that interior degree.

For T, the only single-edge sides are \(DA\) and \(CD\). Their incident
apices have degree five and cannot be flipped. An unflipped joining edge
separates the core graphs on its two sides: no interior core edge crosses
it or uses its boundary endpoints. Thus a later D or a flipped join beyond
an R strip would disconnect the final core graph. Only unflipped R pieces
can follow T; axis locking excludes simultaneous strips on its two
perpendicular eligible sides. The rectangle theorem identifies the old
paths with complete opposite sides. This proves necessity.

Here is a common-weight construction for every listed presentation.
Fix a transverse width \(W>0\). On a maximal run of flipped joins choose
one \(d>0\). Give each participating R height \(2d\). For D set
\[
 h=\frac1{10},\qquad H=\frac d h,\qquad
 k=\frac{1-4h}{4(1-h)},
\]
use the complete positive D chart \eqref{eq:Dcoordinates}.
Both end centers are distance \(hH=d\) from their ends. Thus each designated
join has two midpoint apices at opposite distances \(d\), and its flip
is a genuine centered axial diamond with positive \(E=0\) faces.

Different maximal flipped runs have no piece in common. Across an
unflipped join only the side lengths must agree; the facing apex distances
need not agree. Hence the runs and remaining pieces fit in one rectangle.
An unflipped D using its other end pair is realized by interchanging its
axes, with the appropriate conjugate chart. For T use the complete
positive family \eqref{eq:Tcoordinates} and attach the required unflipped
R strip on its eligible side. Positivity and every face equation of both
families were verified in Proposition~\ref{prop:terminalbases}.

Recover weights on every face from \eqref{eq:weight} in the one final
chart. The shared edge increments give the same midpoint equations;
no link weights are chosen independently. All interior degrees are four
or five, including the apices raised from four to five. Proposition~
\ref{prop:rectangle} gives all vertex gates, both the killed gates at old
boundary vertices and the full-link gates at new vertices. The opposite
old-side marking uses the correspondingly rotated conjugate chart and
newly recovered weights, not an unchanged weighting with a rotated mask.
\end{proof}

\section{Common weights and the original empty-face condition}

Throughout this section \(N=C^0\smallsetminus L^0\) is fixed. Every retained
mass in \(H_N\) is the full sum of the incident face weights. In particular,
new boundary vertices are never removed by an auxiliary marking of a cut.

\begin{lemma}[Positive responses at old interfaces]\label{lem:positivecaps}
Suppose a whole single-edge side \(UV\) of a positive rectangle has an
interior four-valent apex \(s\), and all cap-boundary vertices are old.
There are positive R and D extensions whose legal forward flip preserves
positive definiteness of the original empty matrix. At interface weight
\(\tau\), their Schur responses at \(s\) are respectively
\[
 \delta_R=\frac{29}{13}\tau,\qquad
 \delta_D=\frac{16}{5}\tau.
\]
\end{lemma}
\begin{proof}
Write \(U=(0,0),V=(W,0),s=(W/2,d)\), where \(W,d>0\).
The cap apex must be \(r=(W/2,-d)\). Its recovered interface weight is
\(\tau=2d/(\sqrt3 W)\). Use an R cap of height \(2d\), or a D cap
of height \(10d\) and parameter \(h=1/10\) as above. Thus the scale is
fixed by the actual remainder, rather than independently adjustable.

For D, substitution in \eqref{eq:weight}, with ring-face scale \(\lambda\),
gives
\[
 \begin{gathered}
 t_{ABz}=t_{VUr}=\lambda w,\quad
 t_{BVq}=t_{UAp}=\lambda P,\quad
 t_{zqp}=t_{pqr}=\lambda S,\\
 (P,w,S)=\left(\frac3{1-4h},\frac h{1-h},
                    \frac{1-2h}{1+2h}\right).
 \end{gathered}
\]
The four other faces have weight \(\lambda\). At \(h=1/10\) the
interface condition forces \(\lambda=9\tau\), and the cap matrix on
\((p,q,z,r)\) is
\begin{equation}\label{eq:Dcapmatrix}
 \tau\begin{pmatrix}
 75&-12&-15&-15\\-12&75&-15&-15\\
 -15&-15&25&0\\-15&-15&0&25
 \end{pmatrix}.
\end{equation}
Eliminate the orthogonal core difference, then \(z\), then the core sum.
The pivots are \(87\tau,25\tau,45\tau\), leaving \(15\tau\) on \(r\).
For R the corresponding cap scalar is \(8\tau\).

The old faces \(UVr,UVs\), each of weight \(\tau\), are replaced by
\(rsU,rsV\), each of weight \(3\tau\). Because \(U,V\) are old,
the exact change on the retained coordinates \((r,s)\) is
\[
 \tau\begin{pmatrix}5&-6\\-6&5\end{pmatrix}.
\]
If \(H_{\mathrm{rem}}\) includes its face \(UVs\), the remaining block
is therefore
\[
 \begin{pmatrix}
 d_0+5\tau&-6\tau e_s^{\mathsf T}\\
 -6\tau e_s&H_{\mathrm{rem}}+5\tau e_se_s^{\mathsf T}
 \end{pmatrix},\qquad
 d_0=8\tau\text{ or }15\tau.
\]
Completing the square gives the positive pivot \(d_0+5\tau\) and
\(H_{\mathrm{rem}}+\delta e_se_s^{\mathsf T}\), where
\(\delta=5\tau-36\tau^2/(d_0+5\tau)\). This proves both stated
values and the positive-definite assertion. All masses and the current
remainder are included; no fixed-weight invariance under inverse flips
is asserted.
\end{proof}

\begin{lemma}[Strict seed weights]\label{lem:strictseeds}
Each of R, D and T has common positive weights satisfying \(V\) and
\(H_N\succ0\) for either opposite-side marking.
\end{lemma}
\begin{proof}
For R all boundary vertices are old, so its only new vertex has positive
mass. For D with the four-center ends horizontal and vertical sides old,
\eqref{eq:Dcapmatrix} has positive pivots \(87\tau,25\tau,45\tau,15\tau\).
For the complementary marking take \(h=1/10\) and ring scale one.
The reciprocal face groups are \(P=1/5,w=9,S=3/2\). Both center
pivots are \(25/2\); eliminating them gives on the two cores
\[
 \begin{pmatrix}21/5&-4\\-4&21/5\end{pmatrix}\succ0.
\]
This also supplies an unflipped D using the other end pair in a train.

For T use \eqref{eq:Tcoordinates}. With four inner-face weights normalized
to one, face recovery gives the six outer weights
\[
 (3P,P,Q,3Q,3R,R),\qquad
 P=\frac{2-a}{3a},\quad Q=\frac{1+a}{3(1-a)},\quad
 R=\frac{1-2a}{3}.
\]
With vertical old sides, the original new set is \((p,q,r,U)\), and
\[
 H_T=\begin{pmatrix}
 4P+3&-2&-2&-(P+1)\\
 -2&4Q+3&-2&-(Q+1)\\
 -2&-2&4R+3&0\\
 -(P+1)&-(Q+1)&0&P+Q+1
 \end{pmatrix}.
\]
At \(a=1/5\), multiply all weights by ten. The face vector is
\[
 (90,30,5,15,6,2,10,10,10,10).
\]
The leading principal minors of
\(10H_T\) are \(150,7100,173800,2836500\), proving strict positivity.
The coordinate reflection \((x,y)\mapsto(1-y,1-x)\) exchanges the two
markings and the boundary vertices \(U,V\); it transforms the parameter
to \((1-2a)/(2-a)\). It supplies the complementary witness with \(V\)
retained instead of \(U\). All these weights come from the displayed
positive charts, so Proposition~\ref{prop:rectangle} also supplies \(V\).
\end{proof}

\begin{proposition}[Old interfaces admit strict common weights]
\label{prop:oldports}
An admissible marked presentation in which every join endpoint is old
has common positive weights satisfying \(V\) and \(H_N\succ0\).
\end{proposition}
\begin{proof}
For a serial R/D word, construct from one end. At an unflipped join,
the old endpoints give an exact direct sum on the original new set;
use the appropriate strict seed from Lemma~\ref{lem:strictseeds}.
At a flipped join, attach the next R or D by Lemma~\ref{lem:positivecaps}.
Its height matches the facing four-center, and its positive response
preserves the complete current matrix's positive definiteness.
A D has another unused four-center for the next possible flip, whereas
an R cannot be used twice by admissibility. This constructs one positive
chart and one common weighting for the entire finite word. No equality
of the two face weights is needed at an unflipped join.

For T with an R strip, old join endpoints mean a free-side attachment.
The original matrix is exactly the strict T matrix, including its new
boundary vertex, in direct sum with the positive R scalars. The empty
strip is Lemma~\ref{lem:strictseeds}. The degree bound and the common
chart ensure all vertex conditions throughout.
\end{proof}

\begin{lemma}[The opposite marking is obstructed]\label{lem:newports}
Suppose a positive midpoint realization has an actual admissible serial
R/D decomposition with at least two pieces and with its short outer ends
old. Its original empty matrix is not positive semidefinite. The same
conclusion holds for T with a nonempty R strip in the old-side direction.
\end{lemma}
\begin{proof}
A triangle with \(j\) original new vertices contributes
\(j(2-j)t\) to \(\boldsymbol1_N^{\mathsf T}H_N\boldsymbol1_N\).
Work first with the face ledger before restoring the designated flips,
but keep the final original set \(N\). In module \(i\), let
\(\tau_i>0\) be the common weight on its two end faces. In an R piece
the side-face weights are \(\tau_i/3\). For D denote by
\(\nu_i,\lambda_i,\sigma_i>0\) its other-side, ring and inner weights.
Both permitted D end pairs have these opposite-face equalities; the
argument uses only positivity, so it does not fix the wrong orientation.
Direct face counting gives
\[
\begin{array}{c|cc}
 &R&D\\ \hline
 \text{one global old end}
 &-2\tau_i&-2\tau_i-6\lambda_i-6\sigma_i\\
 \text{no old end}
 &-8\tau_i&-6\tau_i-6\nu_i-12\lambda_i-6\sigma_i.
\end{array}
\]
For example, an end D has one end contribution \(\tau_i\), another
\(-3\tau_i\), two all-new ring faces and two all-new inner faces;
its other four faces have two new vertices and contribute zero.
An internal piece has all its vertices new.

At a genuine centered flip the two equal weights \(\tau\) become
\(\tau/3\) in this marking. All four flip vertices are new, so the
indicator increases by
\[
 -6(\tau/3)+6\tau=4\tau.
\]
Allocate \(2\tau\) to each adjacent piece. An end piece receives at
most one credit, and an internal piece at most two. Hence an end R has
corrected contribution at most zero, an end D strictly negative, and
every internal piece strictly negative. This proves the claim for at
least three pieces, and for two pieces if either is D. For two unflipped
R pieces the total \(-2\tau_1-2\tau_2\) is negative. Equality of
neighboring scales was used only at an actual flipped join, not along
an arbitrary chain of unflipped joins.

The only remaining case is two R pieces with their joining edge flipped.
This is F. Write \(\tau=3r\). On its four original new vertices
\((p,q,U,V)\), the masses are \((7r,7r,3r,3r)\) and every retained edge
has coefficient \(2r\). Thus
\[
 \boldsymbol1^{\mathsf T}H_N\boldsymbol1=0,
 \qquad H_N\boldsymbol1=r(1,1,-1,-1)^{\mathsf T}\ne0,
\]
which is impossible for a positive semidefinite matrix.

Finally take T with a nonempty strip on \(DA\) and vertical old sides.
Its original new vertices now include \(A,D,U\) and the three cores.
Using the six outer weights in Lemma~\ref{lem:strictseeds}, its
contribution to the same global indicator is
\[
 -12P+3Q+3R-9<-17,
 \qquad P>1,\quad 0<Q<1,\quad 0<R<1/3.
\]
Indeed \(DAp,AUp\) contribute \(-9P,-3P\); \(BVq,VCr\) contribute
\(3Q,3R\); three inner faces are all-new; all other contributions
vanish. Each R piece contributes strictly negatively by the same table.
The reflected case changes both the coordinates and the original mask.
This proves the T assertion for all positive parameters and scales.
\end{proof}

\begin{proof}[Completion of the proof of Theorem~\ref{thm:main}]
Let an arbitrary positive complete weighting be given. Theorem~
\ref{thm:reduction} and Proposition~\ref{prop:presentations} supply a
legal serialization of its own positive geometry. In that serialization
the old sides are a complete opposite pair. If there is a join, either
all its endpoints are old, or the realization has the opposite marking
excluded by Lemma~\ref{lem:newports}. Thus the given weighting supplies
an admissible combinatorial presentation whose join endpoints are all old.
The presentation is thus derived from the given weighting itself.

Conversely, from any such combinatorial presentation, Proposition~
\ref{prop:oldports} constructs common positive weights satisfying \(V\)
and \(H_N\succ0\). With no joins the condition is vacuous, and the
three strict seeds apply. For F, the original endpoints of its two-piece
join enter the marking test. Since strict positivity implies semidefiniteness,
this proves the equivalence of the combinatorial criterion, existence
with \(H_N\succeq0\), and existence with \(H_N\succ0\).
This establishes all three equivalent existence statements and completes
the classification.
\end{proof}

\section{Discussion}
Theorem~\ref{thm:main} separates vertex compatibility, the original
empty-face condition, and existence of common positive weights. The
first has a geometric reduction; the second depends on the original
interface marking; the third is constructive in the stated class.
The key structural step is the realization of the core eliminations
by positive geometric reductions.

\paragraph{Disconnected cores: a sufficient family.}
The same positive-response construction also produces strict common
weights for serial supports with disconnected cores.
\begin{corollary}[Strict weights for disconnected serial cores]
\label{cor:disconnectedserial}
Every admissible serial R/D presentation whose original join endpoints
are old admits common positive weights satisfying $V$ and $H_N\succ0$,
without a connectivity assumption on $G_5$.
For any integers $m\ge1$ and $k_1,\ldots,k_m\ge1$, such a presentation
exists with precisely $m$ core components of orders $4k_i-2$.
\end{corollary}
\begin{proof}
The construction in Proposition~\ref{prop:oldports} uses positive cap
responses at flipped joins and direct sums at unflipped old joins;
it never assumes that $G_5$ is connected. For the second assertion,
take runs of $k_i$ D pieces, using $AB$ and $UV$ as their end edges.
Flip every join within each run, and leave the joins between runs
unflipped. Each D contributes its adjacent cores $p,q$; each internal
flip adds two formerly four-valent apices and connects the neighboring
cores. Distinct flips use distinct apices. A run therefore has
$2k_i+2(k_i-1)=4k_i-2$ cores in one component. At a join between runs,
only the old boundary endpoints are identified, and no core edge joins
the two runs. The common-weight construction proves the assertion.
\end{proof}

This gives arbitrarily many arbitrarily large components, but only for
the specified presentations. Necessity for arbitrary disconnected cores
remains open. Connectedness was used to exclude a four-star whose
opposite cores lie on different sides of its boundary-to-boundary
diagonal: a connecting core path would have to cross that diagonal.
With disconnected cores there need be no such path. Likewise an
exterior region cut off by a boundary chord may contain another core
component, so the empty-cap and physical-boundary conclusions of the
reduction argument no longer follow. The missing step is control of
these interfaces while preserving the original new-vertex set and a
single positive face weighting, rather than merely a bound on each
component separately. The interior bound of
Proposition~\ref{prop:blockbudget} itself does not need connectedness.

\paragraph{Other degrees and parameters.}
For higher interior degrees, midpoint equations alone need not imply
all vertex spectral inequalities. At a parameter other than $1/2$,
the weighted identity~\eqref{eq:garland} does not force equality in
each local quadratic bound, so the conjugate-rectangle reduction is
not automatically available. Even in the present degree range,
failure of the support criterion proves nonattainment at $1/2$, not
a uniform gap for the infimum over degenerating positive weights.
These questions, and the passage from relative killed operators to
the absolute operators of Vinzant's question, remain unresolved.

\end{document}